\documentclass[preprint,11pt]{elsarticle}
\usepackage{amssymb}
\usepackage{amsfonts}
\usepackage{amsmath}
\usepackage{enumerate}
\usepackage[dvipdf]{rotating}
\usepackage{rotating}
\usepackage[T1]{fontenc} 
\usepackage{graphicx,scrpage2} 
\usepackage{float,colortbl,xcolor,booktabs,layout} 
\usepackage{hyperref}
\usepackage{psfrag}
\usepackage{epstopdf}

\newtheorem{theorem}{Theorem}

\newtheorem{example}{Example}

\newtheorem{proposition}{Proposition}
\newtheorem{remark}{Remark}

\newenvironment{proof}[1][Proof]{\noindent\textbf{#1.} }{\ \rule{0.5em}{0.5em}\medskip}

\newcommand{\eps}{\varepsilon}
\newcommand{\R}{\mathbb{R}}
\newcommand{\N}{\mathbb{N}}

\newcommand\Fc{\mathcal{F}}

\newcommand\BB{\mathbb{B}}

\newcommand\EE{\mathbb{E}}
\newcommand\GG{\mathbb{G}}
\newcommand\HH{\mathbb{H}}
\newcommand\PP{\mathbb{P}}
\newcommand\Sb{\mathbb{S}}

\newcommand\weak{\ \overset{w.}{{\longrightarrow}} \ }
\newcommand\pconv{\ \overset{\PP}{{\longrightarrow}} \ }
\newcommand\weakcondp[1]{\overset{\PP}{\underset{ #1}{\longrightarrow}}}

\DeclareMathOperator{\Cov}{Cov}
\DeclareMathOperator{\Var}{Var}

\begin{document}

\thispagestyle{empty}
\pagenumbering{Roman}
\setcounter{page}{0}
\renewcommand{\thefootnote}{\fnsymbol{footnote}}
\setcounter{footnote}{1}
\parindent0cm
\renewcommand{\thefootnote}{\fnsymbol{footnote}}

\begin{center}
{\Large {\bf  Consistent testing for a constant copula under strong mixing based on the tapered block multiplier technique}}\\

\bigskip

\bigskip

{\large  Axel B\"ucher\footnote[2]{\label{foot:2}Ruhr-Universit\"at Bochum,
Universit\"atsstra\ss e 150, 44780 Bochum, Germany;  Email: axel.buecher@rub.de.de, \\ Tel: +49 (0) 234 3223286.} and Martin Ruppert\footnote[1]{\label{foot:1}Graduate School of Risk Management and Department of Economic and Social Statistics, University
of Cologne, Albertus-Magnus-Platz, 50923 K\"{o}ln, Germany;  Email: martin.ruppert@uni-koeln.de, \\ Tel: +49 (0) 221 4706656, Fax: +49 (0) 221 4705074.}}

\bigskip

\bigskip

\thispagestyle{empty}

\bigskip

\bigskip

\vspace*{1cm}

{\Large \textbf{Abstract}}
\end{center}
\bigskip
Considering multivariate strongly mixing time series, nonparametric tests for a constant copula with specified or unspecified change point (candidate) are derived; the tests are consistent against general alternatives. A tapered block multiplier technique based on serially dependent multiplier random variables is provided to estimate p-values of the test statistics. Size and power of the tests in finite samples are evaluated with Monte Carlo simulations.  The block multiplier technique might have several other applications for statistical inference on copulas of serially dependent data.

\bigskip

\bigskip


\bigskip



\newpage
\pagenumbering{arabic}
\setcounter{page}{1}
\renewcommand{\thefootnote}{\arabic{footnote}}
\setcounter{footnote}{0}

\section{Introduction}

Over the last decade, copulas have become a standard tool in modern risk management. The copula of a continuous random vector is a function which uniquely determines the  dependence structure linking the marginal distribution functions. Copulas play a pivotal role for, e.g., measuring multivariate association \citep[see][]{schmidetal10}, pricing multivariate options \citep[see][]{vandengoorberghetal05} and allocating financial assets \citep[see][]{patton04}. The latter two references emphasize that time variation of copulas possesses an important impact on financial engineering applications.

Evidence for time-varying dependence structures can indirectly be drawn from functionals of the copula, e.g., Spearman's $\rho$, as suggested by \citet{gaisseretal09} and \citet{wiedetal11}. Investigating time variation of the copula itself, \citet{busettiharvey10} consider a nonparametric quantile-based test for a pointwise constant copula. Semiparametric tests for time variation of the parameter within a prespecified parametric copula model are proposed by \citet{diasembrechts09} and \citet{giacominietal09}. \citet{gueganzhang10} combine tests for constancy of the copula (on a given set of vectors on its domain), the copula family, and the parameter. All the latter references on time variation of the copula are based on the assumption of i.i.d. observations.
With respect to financial time-series, this assumption may be approximated by
the estimation of a GARCH model and by using the residuals obtained after GARCH filtration. However, the effect of replacing unobserved innovations by estimated residuals has to be taken into account. Therefore, specific techniques for residuals are required (cf. \cite{chenfan06}) and exploring this approach, \citet{remillard10} investigates a nonparametric change point test for the copula of residuals in stochastic volatility models.

In the present paper we go beyond these approaches and work in a purely nonparametric setting and under the allowance of general serial dependence of a multivariate time series measured by their alpha-mixing coefficients. In similar settings, \citet{fermanianscaillet03}, \citet{doukhanetal09} and \citet{buecvolg2011} consider the nonparametric estimation of copulas. These references form the basis for our new tests for a constant copula under strong mixing mixing assumptions. We introduce nonparametric Cram\'{e}r-von Mises-, Kuiper-, and Kolmogorov-Smirnov tests which assess constancy of the copula on its entire domain. In consequence, they are consistent under general alternatives. Depending on the object of investigation, tests with a specified or unspecified change point (candidate) are introduced. Whereas the former setting requires a hypothesis on the change point location, it allows us to relax the assumption of strictly stationary univariate processes. P-values of the tests are estimated based on a generalization of the multiplier bootstrap technique introduced in \citet{remillardscaillet09} to the case of strongly mixing time series. This new technique can be used to generalize several statistical inference methods from the i.i.d. case to the serial dependent case.
The idea of our approach is comparable to block bootstrap methods with the difference that, instead of sampling blocks with replacement, we generate blocks of serially dependent multiplier random variables. For a general introduction to the latter idea, we refer to \citet{buehlmann93} and \citet{paparoditispolitis01}.
Independently from our work, a current working paper of \citet{vankampenwied12} investigates similar Kolmogorov-Smirnov tests for constancy of the copula with unspecified change point based on a bootstrap technique introduced in \citet{inoue01}.

This paper is organized as follows. In Section \ref{section:inference}, we discuss weak convergence of the empirical copula process under strong mixing.
We derive a tapered block multiplier bootstrap technique for inference on the weak limit and assess it in finite samples. Tests for a constant copula with specified or unspecified change point (candidate) which are relying on this technique are established in Section \ref{section:tests}. Section \ref{section:conclusion} concludes the paper. All the proofs are deferred to \ref{appendix:proofTheorems}.

\section{Nonparametric inference based on serially dependent observations}
\label{section:inference}

As a basis for the tests introduced in Section \ref{section:tests}, we briefly recapitulate the results of \citet{doukhanetal09}, \citet{segers11} and \citet{buecvolg2011} on the asymptotic behavior of the empirical copula process under nonrestrictive smoothness assumptions in the case of strongly mixing observations.
The main result of this section is
a multiplier-based resampling method for this particular setting. We establish its asymptotic behavior and investigate the performance in finite samples.

\subsection{Asymptotic theory}

Consider a vector-valued {\it process}\index{process} $(\mathbf{X}_j)_{j \in \mathbb{Z}}$ with $\mathbf{X}_j=(X_{j,1}, \ldots, X_{j,d})$ taking values in $\mathbb{R}^d.$ Let $F_{i}$ be the distribution function of $X_{j,i}$ for all $j \in \mathbb{Z},$ $i=1, \ldots, d$ and let $F$ be the joint distribution of $\mathbf{X}_j$ for all $j \in \mathbb{Z}.$ Assume that all marginal distribution functions are continuous. Then, according to Sklar's Theorem \citep{sklar59}, there exists a unique copula $C$ such that $F(x_1, \ldots, x_d)=C(F_1(x_1), \ldots, F_d(x_d))$ for all $(x_1, \ldots, x_d) \in \R^d.$ The $\sigma$-fields generated by $\mathbf{X}_j,j\leq t,$ and $\mathbf{X}_j,j\geq t,$ are denoted by $\mathcal{F}_t=\sigma\{\mathbf{X}_j,j\leq t\}$ and $\mathcal{F}^{t}=\sigma\{\mathbf{X}_j,j\geq t\}$, respectively. We define
\begin{align*}
\alpha(\mathcal{F}_s,\mathcal{F}^{s+r})=\sup_{A\in \mathcal{F}_s,B \in \mathcal{F}^{s+r}}|P(A\cap B)-P(A)P(B)|.
\end{align*}
The strong- (or $\alpha$-) mixing coefficient $\alpha_{\mathbf{X}}$ corresponding to the process $(\mathbf{X}_j)_{j\in \mathbb{Z}}$ is given by
$\alpha_{\mathbf{X}}(r)=\sup_{s \geq 0}\alpha(\mathcal{F}_s,\mathcal{F}^{s+r}).$ The process  $(\mathbf{X}_j)_{j\in \mathbb{Z}}$ is said to be {\it strongly mixing} if $\alpha_{\mathbf{X}}(r)\rightarrow 0$ for $r\rightarrow \infty.$ This type of weak dependence covers a broad range of time-series models. Consider the following examples, cf. \citet{doukhan94} and \citet{carrascochen02}:

\begin{example}
\label{example:mixing}
i) {\it $\mathrm{AR(1)}$ processes} $(X_j)_{j \in \mathbb{Z}}$ given by
\begin{align*}
X_j=\beta X_{j-1}+ \eps_j,
\end{align*}
where $(\eps_j)_{j \in \mathbb{Z}}$ is a sequence of independent and identically distributed continuous innovations with mean zero. For $|\beta|<1$, the process is strictly stationary and strongly mixing with exponential decay of $\alpha_{\mathbf{X}}(r)$.

ii) {\it $\mathrm{GARCH(1,1)}$ processes} $(X_j)_{j \in \mathbb{Z}}$,
\begin{align}
\label{equation:GARCH}
X_j=\sigma_j \eps_j, \quad \sigma_j^2=\omega+\beta \sigma_{j-1}^2+\alpha \eps_{j-1}^2,
\end{align}
where $(\eps_j)_{j \in \mathbb{Z}}$ is a sequence of independent and identically distributed continuous innovations, independent of $\sigma_0^2,$ with mean zero and variance one. For $\alpha+\beta<1$, the process is strictly stationary and strongly mixing with exponential decay of $\alpha_{\mathbf{X}}(r)$.

iii) For multivariate analogues of i) and ii) see Section \ref{subsection:NI_fs} below.
\end{example}

Let $\mathbf{X}_{1},\ldots,\mathbf{X}_{n}$ denote a sample from $(\mathbf{X}_j)_{j \in \mathbb{Z}}.$ A simple nonparametric estimator for the copula $C$ is given by the {\it empirical copula} which is first considered by \citet{rueschendorf76} and \citet{deheuvels79}. Depending on whether the marginal distribution functions are assumed to be known or unknown, we define, for $\mathbf{u}\in \left[ 0,1\right] ^{d}$,
\begin{align}
\label{equation:empirical_copula}
& C_{n}(\mathbf{u}) := \frac{1}{n}\sum_{j=1}^n \mathbf{1}_{\{
\mathbf{U}_{j}\leq \mathbf{u}\}}, &
& \widehat{C}_{n}(\mathbf{u}) := \frac{1}{n}\sum_{j=1}^n \mathbf{1}_{\{
\widehat{\mathbf{U}}_{j}\leq \mathbf{u}\}}, &
\end{align}
where $\mathbf{U}_j=(U_{j,1},\dots,U_{j,d})$ and $\widehat{\mathbf{U}}_j=(\widehat{U}_{j,1},\dots,\widehat{U}_{j,d})$ with observations $U_{j,i}=F_i(X_{j,i})$ and {\it pseudo-observations} $\widehat{U}_{j,i}=\widehat{F}_{i}(X_{j,i})$ for all
$j=1,\dots,n$ and $i=1,\dots,d$, and where $\widehat{F}_{i}(x)=n^{-1}\sum_{j=1}^n\mathbf{1}_{\{X_{j,i}\leq x\}}$
for all $x\in \R.$ Unless otherwise noted, the marginal distribution functions are assumed to be unknown and the {\it empirical copula}  $\widehat{C}_n$ is used. In addition to the practical relevance of this assumption, \citet{genestsegers10} prove that pseudo-observations $\widehat{\mathbf{U}}_{j}$ permit more efficient inference on the copula than observations $\mathbf{U}_{j}$ for a broad class of copulas.

\citet{doukhanetal09} investigate dependent observations and establish the asymptotic behavior of the {\it empirical copula process}, defined by $\sqrt{n}\{\widehat{C}_n-C\},$ assuming the copula to possess continuous partial derivatives on $[0,1]^d.$ \citet{segers11} points out that many popular families of copulas (e.g., the Gaussian, Clayton, and Gumbel--Hougaard families) do not satisfy the assumption of continuous first partial derivatives on $[0,1]^d.$ He establishes the asymptotic behavior of the empirical copula for serially independent observations under the weaker condition
\begin{align}
\label{equation:segers}
D_iC(\mathbf{u}) \ \text{exists and is continuous on} \ \left\{ \mathbf{u} \in [0,1]^d \left| u_i \in (0,1)\right.\right\} \ \text{for all} \ i=1,\ldots, d.
\end{align}
Under Condition (\ref{equation:segers}), the partial derivatives' domain can be extended to $\mathbf{u} \in [0,1]^d$ by
\begin{align}
\label{equation:pd}
D_{i}C(\mathbf{u})= \left\{
\begin{array}{lll} 
\lim_{h \rightarrow 0} \frac{C(\mathbf{u}+h\mathbf{e}_i)-C(\mathbf{u})}{h} & \text{for all} \ u \in [0,1]^d, \ 0<u_i<1, \\
\limsup_{h \downarrow 0} \frac{C(\mathbf{u}+h\mathbf{e}_i)}{h} & \text{for all} \ u \in [0,1]^d, \ u_i=0, \\ \limsup_{h \downarrow 0} \frac{C(\mathbf{u})-C(\mathbf{u}-h\mathbf{e}_i)}{h} & \text{for all} \ u \in [0,1]^d, \ u_i=1, \
\end{array} \right.
\end{align}
and for all $i=1, \ldots, d,$ where $\mathbf{e}_i$ denotes the $i$th column of a $d \times d$ identity matrix. The following generalization of the results in \cite{doukhanetal09} and \cite{segers11} is a consequence of Theorem 2.4 in \citet{buecvolg2011}, see Corollary 2.5 in that reference.

\begin{theorem} \label{theorem:empcopprocess}
Consider observations $\mathbf{X}_{1},\ldots,\mathbf{X}_{n},$ drawn from a strictly stationary process $(\mathbf{X}_j)_{j \in \mathbb{Z}}$ satisfying the strong mixing condition $\alpha_\mathbf{X}(r)=\mathcal{O}(r^{-a})$ for some $a>1$. Then
\begin{align*}	
	\BB_{C,n} = \sqrt{n}(C_n-C) \weak \BB_C
\end{align*}
in the metric space space of uniformly bounded functions on $[0,1]^d$ equipped with the uniform metric $\left(\ell^\infty([0,1]^{d}),\|\cdot\|_\infty \right)$. Here, $\BB_C$ denotes a centered tight Gaussian process on $[0,1]^d$ with covariance function
\begin{equation}
\label{equation:covariancestructure}
	\gamma(\mathbf{u},\mathbf{v}) = \Cov(\BB_{C}(\mathbf{u}),\BB_{C}(\mathbf{v}))=\sum_{j \in \mathbb{Z}} \Cov\left(\mathbf{1}_{\{\mathbf{U}_0 \leq \mathbf{u}\}},\mathbf{1}_{\{\mathbf{U}_j \leq \mathbf{v} \}}\right) \ \mbox{for all} \ \mathbf{u},\mathbf{v} \in [0,1]^d.
\end{equation}
Moreover, if $C$ satisfies Condition (\ref{equation:segers}), then
\begin{align*}
	\widehat{\GG}_{C,n}=\sqrt{n}(\widehat{C}_n -C) \weak \GG_C
\end{align*}
in $\ell^\infty([0,1]^{d})$, where $\GG_C$ represents a  Gaussian process given by
\begin{equation}
\label{equation:limitingprocess}
	\GG_{C}(\mathbf{u})=\BB_{C}(\mathbf{u})- \sum_{i=1}^d D_{i}C(\mathbf{u})\BB_{C}(\mathbf{u}^{(i)}) \ \mbox{for all} \ \mathbf{u} \in [0,1]^d.
\end{equation}
Here, $\mathbf{u}^{(i)}$ denotes the vector where all
coordinates, except the $i$th coordinate of $\mathbf{u}$, are
replaced by $1$.
\end{theorem}
Notice that the covariance structure as given in Equation (\ref{equation:covariancestructure}) depends on the entire process $(\mathbf{X}_j)_{j \in \mathbb{Z}}$ in case it is not serially uncorrelated.

\subsection{Resampling techniques}
\label{subsection:resampling}

In this Section, we introduce two bootstrap techniques for the empirical copula process which are applicable in the case of strongly mixing observations.
We begin with the (moving) block bootstrap, which serves as a benchmark in the finite sample assessment. Subsequently, we derive the main result of this Section about the asymptotic consistency of a generalized multiplier bootstrap technique.


\subsubsection{The block bootstrap}
\citet{fermanianetal04} investigate the empirical copula process for independent and identically distributed observations $\mathbf{X}_1,\ldots, \mathbf{X}_n$ and prove consistency of the {\it nonparametric bootstrap}\index{bootstrap} method which is based on sampling with replacement from $\mathbf{X}_1,\ldots, \mathbf{X}_n.$ We denote a bootstrap sample by $\mathbf{X}^{B}_1, \ldots, \mathbf{X}^{B}_n$ and define
\begin{align*}
& \widehat{C}^{B}_n(\mathbf{u}):=\frac{1}{n} \sum_{j=1}^n \mathbf{1}_{ \left\{ \widehat{\mathbf{U}}^{B}_j \leq \mathbf{u}\right\}} \ \mbox{for all} \ \mathbf{u} \in [0,1]^d,
& \widehat{U}^{B}_{j,i}:=\frac{1}{n} \sum_{k=1}^n \mathbf{1}_{ \left\{ X^{B}_{k,i} \leq X^{B}_{j,i} \right\} }
\end{align*}
for all $j=1,\ldots,n$ and $i=1,\ldots,d.$ Notice that the bootstrap empirical copula can equivalently be expressed based on multinomially $(n, n^{-1}, \ldots, n^{-1})$ distributed random variables $\mathbf{W}=(W_1, \ldots, W_n):$
\begin{align*}
& \widehat{C}^{W}_n(\mathbf{u}):=\frac{1}{n} \sum_{j=1}^n W_j \mathbf{1}_{ \left\{ \widehat{\mathbf{U}}^W_j \leq \mathbf{u} \right\} } \ \mbox{for all} \ \mathbf{u} \in [0,1]^d,
& \widehat{U}^{W}_{j,i}:=\frac{1}{n} \sum_{k=1}^n W_k \mathbf{1}_{ \left\{ X_{k,i} \leq X_{j,i} \right\} }
\end{align*}
for all $j=1,\ldots,n$ and $i=1,\ldots,d.$ It is shown in \cite{fermanianetal04, buecher11} that the corresponding {\it bootstrap empirical copula process} {\it converges weakly conditional on $\mathbf{X}_1, \ldots, \mathbf{X}_n$ in probability} in $\left(\ell^\infty([0,1]^{d}),\|\cdot\|_\infty \right)$, notationally
\begin{equation}
\label{equation:resboot}
	\widehat{\GG}_{C,n}^W =\sqrt{n}(\widehat{C}^W_n-C_n) \weakcondp{W} \GG_C.
\end{equation}
Here, \textit{weak convergence conditional on the data in probability} is understood in the Hoffmann-J\o rgensen sense as defined in \citet{kosorok08}, i.e., $\widehat{\GG}_{C,n}^W  \weakcondp{W} \GG_C$ if and only if the following two conditions hold:
\begin{align}
\label{equation:cwc1}
&\sup_{h \in BL_1(\ell^\infty([0,1]^d))} \left|\mathbb{E}_Wh(\widehat{\GG}_{C,n}^W)-\mathbb{E}h\left(\GG_C\right)\right| \pconv 0, \\
\label{equation:cwc2}
&\mathbb{E}_Wh(\widehat{\GG}_{C,n}^W)^*-\mathbb{E}_Wh(\widehat{\GG}_{C,n}^W)_* \pconv 0 \quad\text{for all } h\in BL_1(\ell^\infty([0,1]^d)).
\end{align}
Here, $\pconv$ denotes convergence in outer probability and $\mathbb{E}_W$ denotes expectation with respect to $\mathbf{W}$ conditional on $\mathbf{X}_1, \ldots, \mathbf{X}_n.$ Furthermore, $BL_1(\ell^\infty([0,1]^d))$ denotes the set
\begin{align*}
	\bigg\{f:\ell^\infty( [0,1]^d)\rightarrow \mathbb{R} \mid \|f\|_\infty \leq 1, |f(\beta)-f(\gamma)|\leq\|\beta-\gamma\|_\infty \ \text{for all} \ \gamma, \beta \in \ell^\infty([0,1]^d)\bigg\}.
\end{align*}
of all Lipschitz-continuous functions bounded by 1 with Lipschitz-constant not exceeding 1 and the asterisks in \eqref{equation:cwc2} denote measurable majorants and minorants with respect to the joint data (i.e., $\mathbf{X}_1, \ldots, \mathbf{X}_n$ and $\mathbf{W}$).
Weak convergence conditional on $\mathbf{X}_1, \ldots, \mathbf{X}_n$ almost surely is defined analogously by replacing outer probability convergence in \eqref{equation:cwc1} and \eqref{equation:cwc2} by outer almost sure convergence. Due to the lack of a general continuous mapping Theorem and an easy functional delta method, see \cite{vaartwellner96, kosorok08}, we do not consider the outer almost sure version in this paper.

Whereas the bootstrap is consistent for i.i.d. samples, consistency generally fails for serially dependent samples.Therefore, a {\it block bootstrap method} is proposed by \citet{kuensch89}. Given the sample $\mathbf{X}_1,\dots,\mathbf{X}_n,$ the block bootstrap method requires blocks of size $l_B=l_B(n)$, $l_B(n)\rightarrow \infty$ as $n\rightarrow \infty$ and $l_B(n)=o(n),$ consisting of consecutive observations
\begin{align*}
B_{h,l_B}=\{\mathbf{X}_{h+1},\dots,\mathbf{X}_{h+l_B}\}, \ \text{for all} \ h=0,\dots, n-l_B.
\end{align*}
We assume $n=kl_B$ (otherwise the last block is truncated) and simulate $\mathbf{H}=(H_1,\dots,H_k)$ independent and uniformly distributed random variables on $\{0,\dots,n-l_B\}.$ The block bootstrap sample is given by the observations of the $k$ blocks $B_{H_1,l_B},\dots,B_{H_k,l_B},$ i.e.,
\begin{align*}
\mathbf{X}_{H_1+1},\dots,\mathbf{X}_{H_1+l_B},\mathbf{X}_{H_2+1},\dots,\mathbf{X}_{H_2+l_B}, \dots \dots, \mathbf{X}_{H_k+1}, \dots, \mathbf{X}_{H_k+l_B}.
\end{align*}
Denote the block bootstrap empirical copula based on this sample by $\widehat{C}_n^{B}(\mathbf{u}).$ Its asymptotic behavior can be established
by means of an asymptotic result on the block bootstrap for general $d$-dimensional distribution functions established by \citet[Theorem 3.1]{buehlmann93}, see also Example 2.10 in \citet{buecvolg2011}.

\begin{theorem}
\label{theorem:blockbootstrap}
Consider observations $\mathbf{X}_1, \ldots, \mathbf{X}_n,$ drawn from a strictly stationary process $(\mathbf{X}_j)_{j \in \mathbb{Z}}$ satisfying $\sum_{r=1}^\infty (r+1)^{16(d+1)} \sqrt{\alpha_\mathbf{X}(r)} < \infty.$ Assume that $l_B(n)=\mathcal{O}(n^{1/2-\eps})$ for $0< \eps <1/2.$
If $C$ satisfies Condition (\ref{equation:segers}), then the block bootstrap empirical copula process converges weakly conditional on $\mathbf{X}_1, \ldots, \mathbf{X}_n$ in probability in $\left(\ell^\infty([0,1]^{d}),\|\cdot\|_\infty \right):$
\begin{align*}
\GG_{C,n}^B(\mathbf{u})=\sqrt{n}\left\{\widehat{C}_n^B(\mathbf{u})-\widehat{C}_n(\mathbf{u})\right\} \weakcondp{H} \GG_C(\mathbf{u}).
\end{align*}
\end{theorem}

A brief remark on the condition on the mixing rate is in order here. The summability conditions forces $\alpha_\mathbf{X}(r)$ to be of order $O(r^{-96-\eps})$ for dimension $d=2$, which is far away from the (sharp) rate $O(r^{-1-\eps})$ needed in Theorem \ref{theorem:empcopprocess}. This discrepancy is due to the fact that the literature does not provide stronger results on the consistency of the block bootstrap for the $d$-variate empirical process except the ones in \cite{buehlmann93}, at least to the best of our knowledge. Exploiting more recent techniques (which are beyond of the scope of this paper) we believe that it is possible to get better rates which are comparable to those of the non-bootstrap version, see also \cite{peligrad1998} and \cite{radulovic1998}. The proof of Theorem \ref{theorem:blockbootstrap} being based on the functional delta method would easily transfer these rates to the rank-based copula setting considered in the present paper. Also note that for applications the message is not too bad: most time series models have exponentially decreasing alpha mixing coefficients, see Example \ref{example:mixing} or more precisely the examples in \cite{buecvolg2011}.


\subsubsection{The multiplier bootstrap}
A process related to the bootstrap empirical copula process defined in \eqref{equation:resboot} can be formulated if both the assumption of multinomially distributed random variables is dropped and the marginal distribution functions are left unaltered during the resampling procedure. The resulting bootstrap scheme is known for the i.i.d. context as the \textit{multiplier bootstrap} or the \textit{multiplier method}, see \cite{scaillet05, remillardscaillet09, buecherdette09, segers11}. In the present section we briefly summarize this concept and extend it to the serially dependent case.

Consider i.i.d.\ multiplier random variables $\xi_1, \ldots, \xi_n$ with mean and variance 1, additionally satisfying $\| \xi_j \|_{2,1}:= \int_0^{\infty} \sqrt{P(|\xi_j| > x)} dx < \infty$ for all $j=1, \ldots, n$ (where the last condition is slightly stronger than that of a finite second moment). Replacing the multinomial multiplier random variables $W_1, \ldots, W_n$ by $\xi_1/\bar{\xi}, \ldots, \xi_n/\bar{\xi}$, where $\bar\xi=n^{-1}\sum_{j=1}^n\xi_j$, (ensuring realizations having arithmetic mean one) yields the {\it multiplier (empirical copula) process} which converges weakly conditional on $\mathbf{X}_1, \ldots, \mathbf{X}_n$ in probability in $\left(\ell^\infty([0,1]^{d}),\|\cdot\|_\infty \right)$ in the i.i.d. situation, see \citet{buecherdette09} and more precisely in \citet[Theorem 2.3]{buedetvol2011}:
\begin{align*}
\widehat{\BB}_{C,n}^\xi(\mathbf{u})=\sqrt{n}\left\{\frac{1}{n}\sum_{j=1}^n \frac{\xi_j}{\bar{\xi}} \mathbf{1}_{ \{ \widehat{\mathbf{U}}_j \leq \mathbf{u} \} } - \widehat{C}_n(\mathbf{u})\right\} \overset{\mathbb{P}}{\underset{\xi}{\longrightarrow}} \BB_C(\mathbf{u}).
\end{align*}
For general considerations of multiplier empirical processes, we refer to the monographs of \citet{vaartwellner96} and \citet{kosorok08}.
The process is introduced by \citet{scaillet05} in a bivariate context, a general multivariate version and its unconditional weak convergence are investigated by \citet{remillardscaillet09} and \citet{segers11}. In order to get approximations of the limiting process $\GG_C$ these authors propose to estimate the partial derivatives of the copula in \eqref{equation:limitingprocess} by some estimator $\widehat{D_iC_n}$ (see, e.g., Example \ref{example:pdest} below) and define
\begin{align*}
	\widehat{\GG}_{C,n}^\xi (\mathbf{u}) = \widehat{\BB}_{C,n}^\xi (\mathbf{u}) - \sum_{i=1}^d \widehat{D_iC_{n}}(\mathbf{u})\widehat{\BB}_{C,n}^\xi\left(\mathbf{u}^{(i)}\right),
\end{align*}
which converges to $\GG_C$ conditionally on the data in probability, see \cite{buecherdette09, buecher11}. Here, the estimator $\widehat{D_iC_n}$ is supposed to satisfy the
following two assumptions, which will also be necessary in the serially dependent situation:

\bf C1 \rm There exists a constant $K$ such that $\|\widehat{D_iC_n}\|_\infty\le K$ for all $n\in\N$.

\bf C2 \rm For all $\delta\in(0,1/2)$ one has
\[
	\sup_{\mathbf{u}\in[0,1]^d:u_i\in[\delta,1-\delta]} \left| \widehat{D_iC_n}(\mathbf{u})- D_iC(\mathbf{u}) \right| \pconv 0.
\]


\begin{example}\label{example:pdest}
It is easily seen that
finite differencing yields a simple nonparametric estimator for the first order partial derivatives $D_iC(\mathbf{u})$ which satisfies \bf C1 \rm and \bf C2\rm. More precisely, we define
\begin{align}
\label{equation:finite_differencing}
\widehat{D_{i}C}(\mathbf{u})= \left\{
\begin{array}{lll}
\frac{\widehat{C}_n(\mathbf{u}+h\mathbf{e}_i)-\widehat{C}_n(\mathbf{u}-h\mathbf{e}_i)}{2h} & \text{for all} \ u \in [0,1]^d, \ h \leq u_i \leq 1-h, \\
\frac{\widehat{C}_n(\mathbf{u}+2h\mathbf{e}_i)}{2h} & \text{for all} \ u \in [0,1]^d, \ 0 \leq u_i < h, \\ \frac{\widehat{C}_n(\mathbf{u})-\widehat{C}_n(\mathbf{u}-2h\mathbf{e}_i)}{2h} & \text{for all} \ u \in [0,1]^d, \ 1-h < u_i \leq 1, \
\end{array} \right.
\end{align}
where $h=h_n\to 0$ such that $\inf_n h_n\sqrt{n}>0$ and where $\mathbf{e}_i$ denotes the $i$th column of the $d \times d$ identity matrix \citep[see][]{segers11,buecher11}.
\end{example}

\citet{buecherdette09} find that the multiplier technique yields more precise results than the nonparametric bootstrap in mean as well as in mean squared error when estimating the asymptotic covariance of the empirical copula process in the i.i.d. context. Motivated by the fact that this technique is inconsistent when applied to serially dependent samples, we derive a generalization of the multiplier technique in the following. \citet{inoue01} develops a block multiplier process for general distribution functions based on dependent data in which the same multiplier random variable is used for a block of observations and the composition of blocks remains unaltered throughout the procedure. We  consider a technique in which the composition of blocks is not fixed. More precisely, a {\it tapered block multiplier (empirical copula) process} is introduced based on the work of \citet[][Chapter 3.3]{buehlmann93} and \citet{paparoditispolitis01}: The main idea is to consider a sample $\xi_{1,n}, \ldots, \xi_{n,n}$ from a process $(\xi_{j,n})_{j \in \mathbb{Z}}$ of serially dependent {\it tapered block multiplier random variables}, satisfying:

\vspace{4pt}

\bf A1 \rm $(\xi_{j,n})_{j \in \mathbb{Z}}$ is independent of the observation process $(\mathbf{X}_j)_{j\in \mathbb{Z}}.$

\bf A2 \rm $(\xi_{j,n})_{j \in \mathbb{Z}}$ is a positive {\it$c \cdot l(n)$-dependent} process,\index{process! $l(n)$-dependent}\index{$l(n)$-dependent} i.e., for fixed $j \in \mathbb{Z},$ $\xi_{j,n}$ is independent of $\xi_{j+h,n}$ for all $|h| \geq c \cdot l(n),$ where $c$ is a constant and $l(n) \rightarrow \infty$ as $n \rightarrow \infty$ while $l(n)=o(n).$

\bf A3 \rm \rm $(\xi_{j,n})_{j \in \mathbb{Z}}$ is strictly stationary. For all $j,h \in \mathbb{Z},$ assume $E[\xi_{j,n}]=\mu>0,$ $\Cov[\xi_{j,n},\xi_{j+h,n}]=\mu^2v(h/l(n))$ and $v$ is a bounded function symmetric about zero; without loss of generality, we consider $\mu=1$ and $v(0)=1.$ All central moments of $\xi_{j,n}$ are supposed to be bounded.

\vspace{4pt}

The following theorem is the main result of this section. Weak convergence of the tapered block multiplier process conditional on a sample $\mathbf{X}_1, \ldots, \mathbf{X}_n$ is established. Regarding the strong assumptions on the mixing rate the remark after Theorem \ref{theorem:blockbootstrap} holds true here as well.

\begin{theorem}
\label{theorem:taperedmultiplier}
Consider observations $\mathbf{X}_1, \ldots, \mathbf{X}_n$ drawn from a strictly stationary process $(\mathbf{X}_j)_{j \in \mathbb{Z}}$ satisfying $\sum_{r=1}^\infty (r+1)^c \sqrt{\alpha_\mathbf{X}(r)} < \infty,$ where $c=max\{8d+12, \lfloor 2 / \eps \rfloor +1\}.$ Let the tapered block multiplier process $(\xi_{j,n})_{j \in \mathbb{Z}}$ satisfy \bf A1\it, \bf A2\it, \bf A3 \it with block length $l(n) \rightarrow \infty,$ where $l(n)=\mathcal{O}(n^{1/2-\eps})$ for $0< \eps <1/2.$ Then,
\begin{align*}
	\widehat{\BB}_{C,n}^M (\mathbf{u}) = \sqrt{n} \left( \frac{1}{n} \sum_{j=1}^{n} \frac{\xi_{j,n}}{\bar{\xi}_n} \mathbf{1}_{\left\{\mathbf{\widehat{U}}_{j} \leq \mathbf{u}\right\}}-\widehat{C}_{n}(\mathbf{u})\right) \weakcondp{\xi} \BB_C(\mathbf{u}) \quad \text{in }\ell^\infty([0,1]^{d}),
\end{align*}
where $\bar \xi_n=n^{-1}\sum_{j=1}^n \xi_{j,n}$. Moreover, if \eqref{equation:segers} holds and if $\widehat{D_iC_{n}}$ satisfies conditions \textbf{C1} and \textbf{C2}, then
\begin{align*}
	\widehat{\GG}_{C,n}^M (\mathbf{u}) = \widehat{\BB}_{C,n}^M (\mathbf{u}) - \sum_{i=1}^d \widehat{D_iC_{n}}(\mathbf{u})\widehat{\BB}_{C,n}^M\left(\mathbf{u}^{(i)}\right) \weakcondp{\xi} \GG_C(\mathbf{u}) \quad \text{in }\ell^\infty([0,1]^{d}).
\end{align*}
\end{theorem}

\begin{remark}
\label{remark:centredmultipliers}
The multiplier random variables can as well be assumed to be centered around zero \citep[cf.][Proof of Theorem 2.6]{kosorok08}. Define $\xi_{j,n}^0:=\xi_{j,n}-\mu.$ Then
\begin{align*}
\widehat{\BB}_{C,n}^{M,0}(\mathbf{u})=\frac{1}{\sqrt{n}}\sum_{j=1}^n \left(\xi_{j,n}^0-\bar{\xi}_n^0\right) \mathbf{1}_{\left\{\mathbf{\widehat{U}}_{j} \leq \mathbf{u}\right\}} = \bar{\xi}_n \frac{1}{\sqrt{n}}\sum_{j=1}^n \left(\frac{\xi_{j,n}}{\bar{\xi}_n}-1 \right) \mathbf{1}_{\left\{\mathbf{\widehat{U}}_{j} \leq \mathbf{u}\right\}}=\bar{\xi} \widehat{\BB}^{M}_{C,n}(\mathbf{u})
\end{align*}
for all $ \mathbf{u} \in [0,1]^d.$ This is an asymptotically equivalent form of the above tapered block multiplier process:
\begin{align*}
\sup_{[0,1]^d} \left| \widehat{\BB}^{M,0}_{C,n}(\mathbf{u})- \widehat{\BB}^{M}_{C,n}(\mathbf{u})\right|=
\sup_{[0,1]^d} \left| \left(\bar{\xi}_n-1\right)\widehat{\BB}^{M}_{C,n}(\mathbf{u})\right| \to 0
\end{align*}
almost surely, since $\widehat\BB^{M}_{C,n}(\mathbf{u})$ tends to a tight centered Gaussian limit, unconditionally. The assumption of centered multiplier random variables is abbreviated as \bf A3b \rm in the following.
\end{remark}

In practice, given observations $\mathbf{X}_1, \ldots, \mathbf{X}_n$ of a strictly stationary process $(\mathbf{X}_j)_{j \in \mathbb{Z}}$ satisfying the assumptions of Theorem \ref{theorem:taperedmultiplier}, approximating $\GG_C$ requires the following three steps:
\begin{enumerate}
	\item Estimate the partial derivatives of $C$ by an estimator satisfying conditions \bf C1 \rm and \bf C2\rm, for instance by the estimator given in Example \ref{example:pdest} above.
	\item For $s=1, \ldots,S$ with $S\in\N$, simulate samples $\xi_{1,n}^{(s)},\ldots, \xi_{n,n}^{(s)}$ from a tapered block multiplier process $(\xi_{j,n})_{j\in\mathbb{Z}}$ satisfying \bf A1\rm, \bf A2\rm, \bf A3\rm. For each $s\in\{1,\dots,S\}$, calculate
\begin{align} \label{equation:tapmult_B}
	\widehat{\BB}_{C,n}^{M(s)}(\mathbf{u})=\frac{1}{\sqrt{n}} \sum_{j=1}^{n} \left(\frac{\xi_{j,n}^{(s)}}{\bar{\xi}_n^{(s)}}-1 \right) \mathbf{1}_{\{\mathbf{\widehat{U}}_{j} \leq \mathbf{u}\}} \ \mbox{for all} \ \mathbf{u} \in [0,1]^d.
\end{align}
	\item For each $s\in\{1,\dots,S\}$, calculate
\begin{align} \label{equation:Ghatmult}
	\widehat{\GG}^{M(s)}_{C,n}(\mathbf{u})=\widehat{\BB}_{C,n}^{M(s)}(\mathbf{u})-\sum_{i=1}^d \widehat{D_iC_{n}}(\mathbf{u})\widehat{\BB}_{C,n}^{M(s)}\left(\mathbf{u}^{(i)}\right) \ \mbox{for all} \ \mathbf{u} \in [0,1]^d.
\end{align}
The sample $\{\widehat{\GG}^{M(s)}_{C,n}\}_{s=1,\dots,S}$ is an approximate sample of $\GG_C$.
\end{enumerate}


There are numerous ways to define tapered block multiplier processes $(\xi_{j,n})_{j\in \mathbb{Z}}$ satisfying the assumptions \bf A1\rm, \bf A2\rm, \bf A3 \rm made in Theorem \ref{theorem:taperedmultiplier}.
In the remaining part of this section, a basic version having uniform weights and a refined version with triangular weights are investigated and compared.

\begin{figure}[t]
    \centering
        \begin{minipage}{.5\linewidth}
            \includegraphics[width=\linewidth]{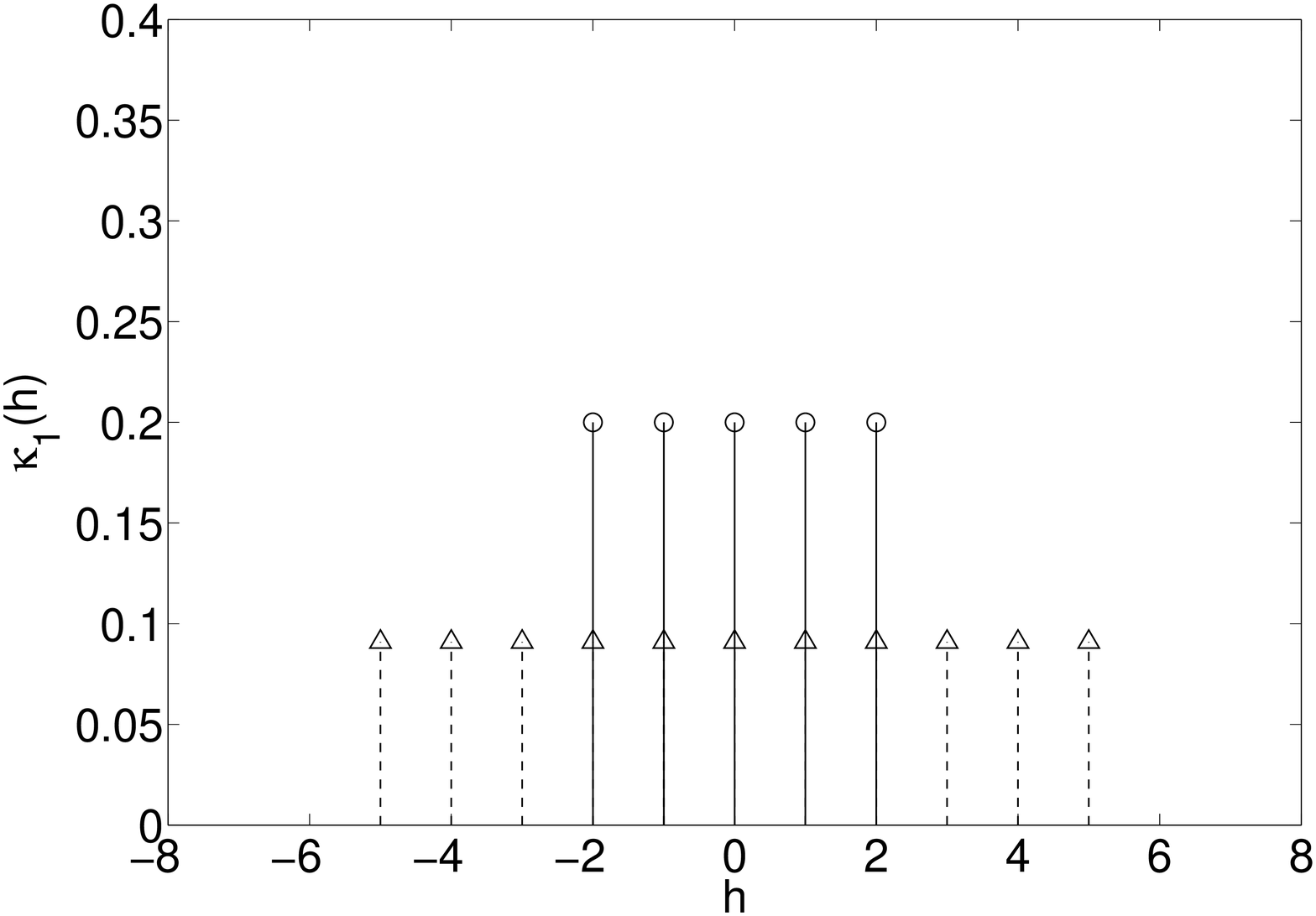}
            \end{minipage}
            \hspace{-.025\linewidth}
        \begin{minipage}{.5\linewidth} 
            \includegraphics[width=\linewidth]{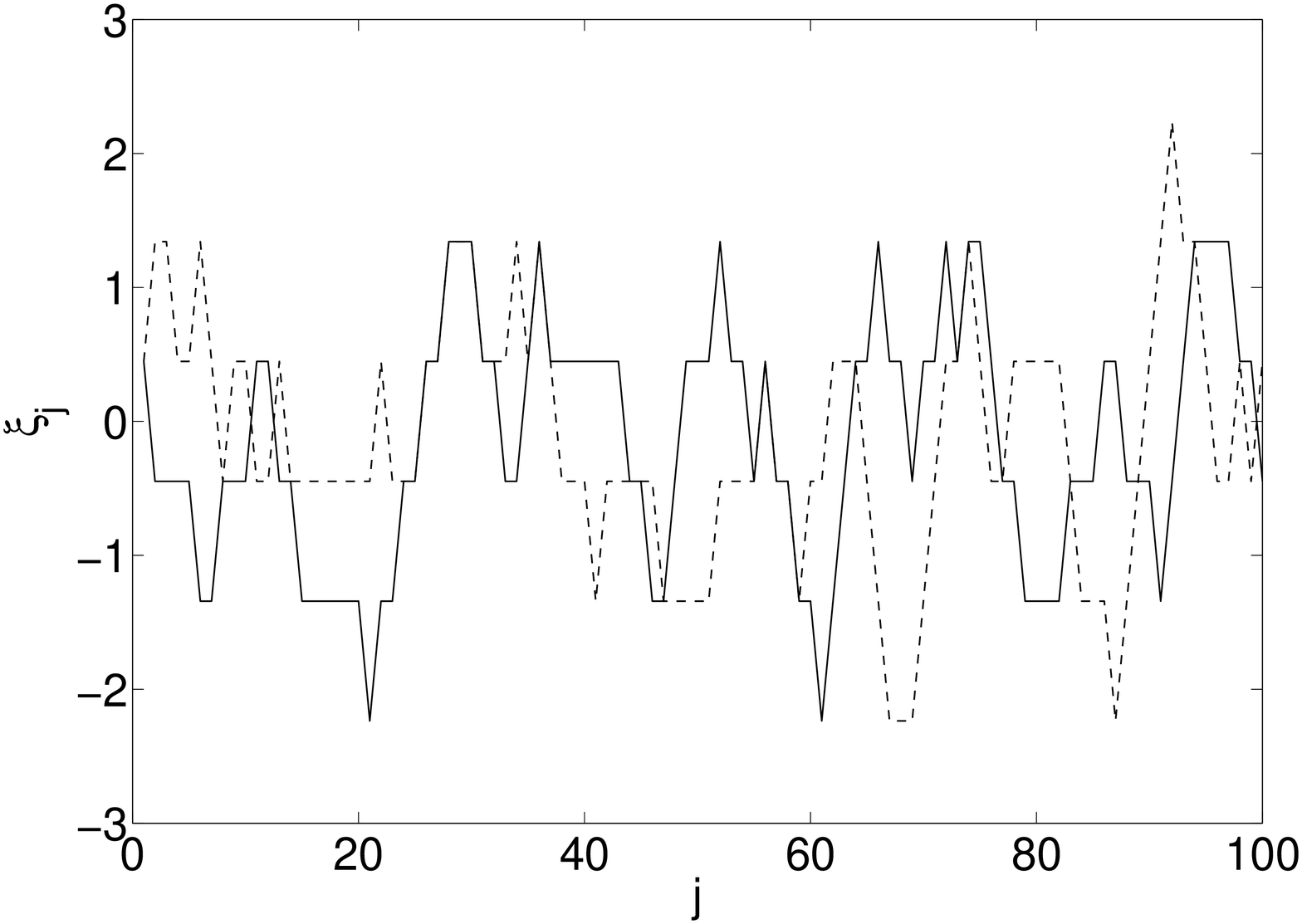}
            \end{minipage}
            \vspace{-10pt}
                \caption{\small{{\bf Tapered block multiplier Monte Carlo simulation}. Kernel function $\kappa_1(h)$ (left) and simulated trajectories of Rademacher-type tapered block multiplier random variables $\xi_1,\ldots, \xi_{100}$ (right) with block length $l(n)=3$ (solid line) and $l(n)=6$ (dashed line), respectively.}}
\label{fig:taperedmultipliers1}
\end{figure}

\begin{example}
\label{example:uniform_blocks}
A simple form of the tapered block multiplier random variables can be defined based on moving average processes. Consider the function $\kappa_1$ which assigns uniform weights given by
\begin{align*}
\kappa_1(h):= \begin{cases}
\frac{1}{2l(n)-1} & \text{for all} \ |h| < l(n) \cr
0 &  \text{else.}
\end{cases}
\end{align*}
Note that $\kappa_1$ is a discrete kernel, i.e., it is symmetric about zero and $\sum_{h\in \mathbb{Z}} \kappa_1(h)=1.$ The tapered block multiplier process is defined by
\begin{align} \label{equation:weightprocess}
	\xi_j=\xi_{j,n}=\sum_{h=-\infty}^\infty \kappa_1(h)w_{j+h} \ \mbox{for all} \ j \in \mathbb{Z},
\end{align}
where $(w_j)_{j\in \mathbb{Z}}$ is an independent and identically distributed sequence of, e.g., Gamma(q,q) random variables with $q:=1/[2l(n)-1].$ The expectation of $\xi_j$ is then given by $E[\xi_j]=1,$ its variance by $Var[\xi_j]=1$ for all $j\in \mathbb{Z}.$ For all $j \in \mathbb{Z}$ and $|h|<2l(n)-1,$ direct calculations further yield the covariance function $\Cov(\xi_j,\xi_{j+h})=\{2l(n)-1-|h|\}/\{2l(n)-1\}$ which linearly decreases as $h$ increases in absolute value. The resulting sequence $(\xi_j)_{j\in \mathbb{Z}}$ satisfies \bf A1\it, \bf A2, \it and \bf A3\it. Exploring Remark \ref{remark:centredmultipliers}, tapered block multiplier random variables can as well be defined based on sequences $(w_j)_{j\in \mathbb{Z}}$ of, e.g., Rademacher-type random variables $w_j$ characterized by $P(w_j=-1/\sqrt{q})=P(w_j=1/\sqrt{q})=0.5$ or Normal random variables $w_j \sim N(0,1/\sqrt{q}).$ In either one of these two cases, the resulting sequence $(\xi_j)_{j\in \mathbb{Z}}$ satisfies \bf A1\it, \bf A2, \it and \bf A3b\it. Figure \ref{fig:taperedmultipliers1} shows the kernel function $\kappa_1$ and simulated trajectories of Rademacher-type tapered block multiplier random variables.
\end{example}

\begin{figure}[ht]
    \centering
        \begin{minipage}{.5\linewidth}
            \includegraphics[width=\linewidth]{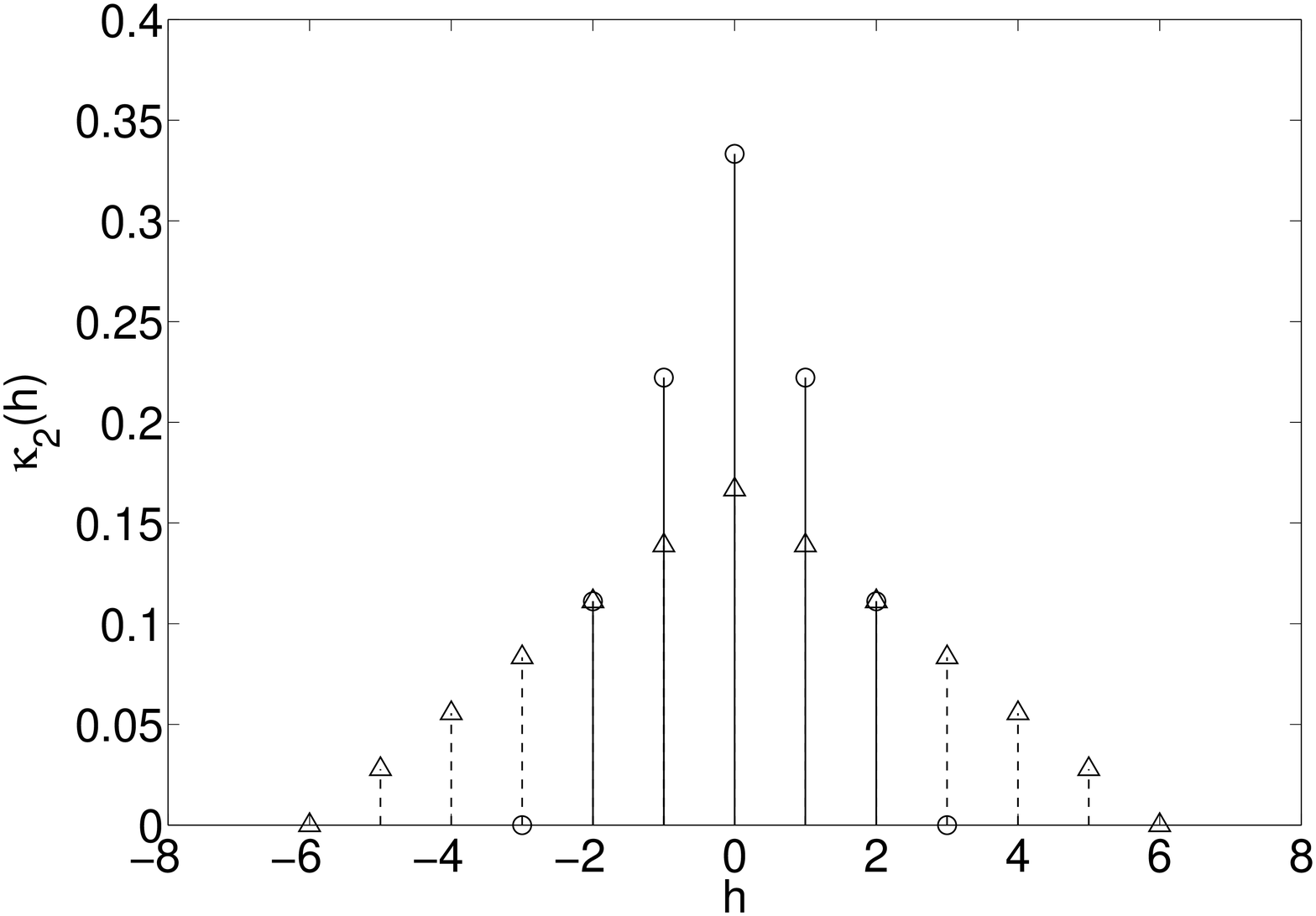}
            \end{minipage}
            \hspace{-.025\linewidth}
        \begin{minipage}{.5\linewidth} 
            \includegraphics[width=\linewidth]{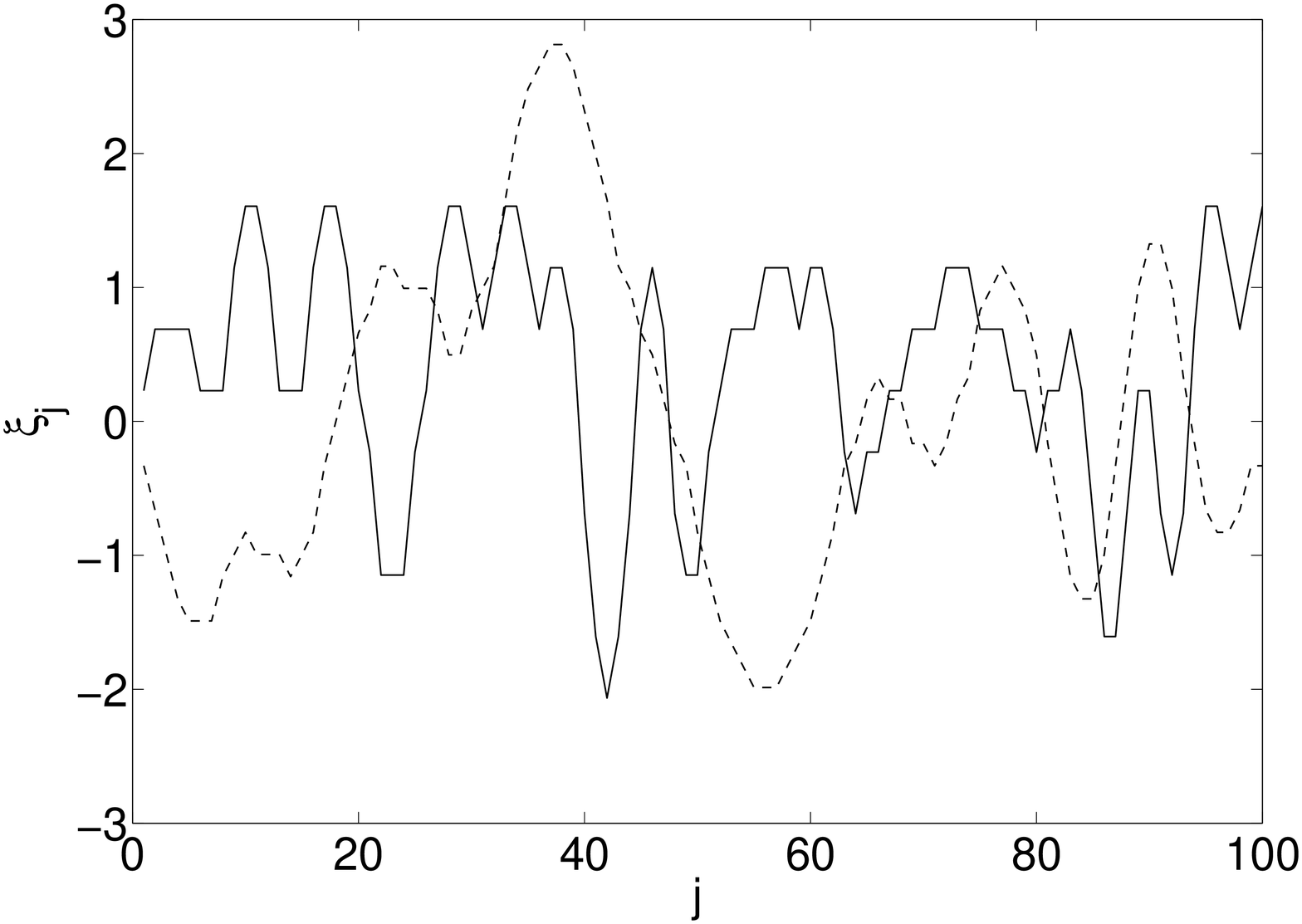}
            \end{minipage}
            \vspace{-10pt}
                \caption{\small{{\bf Tapered block multiplier Monte Carlo simulation}. Kernel function $\kappa_2(h)$ (left) and simulated trajectories of Rademacher-type tapered block multiplier random variables $\xi_1,\ldots, \xi_{100}$ (right) with block length $l(n)=3$ (solid line) and $l(n)=6$ (dashed line), respectively.}}
\label{fig:taperedmultipliers2}
\end{figure}

\begin{example}
\label{example:triangular_blocks}
Following \citet{buehlmann93}, let us define the kernel function by
\begin{align*}
\kappa_2(h):=\max\{0,\{1-|h|/l(n)\}/l(n)\}
\end{align*}
for all $h \in \mathbb{Z}.$ The tapered multiplier process $(\xi_j)_{j\in \mathbb{Z}}$ follows Equation (\ref{equation:weightprocess}), where $(w_j)_{j\in \mathbb{Z}}$ is an independent and identically distributed sequence of Gamma(q,q) random variables with $q=2/\{3l(n)\}+1/\{3l(n)^3\}.$ The expectation of $\xi_j$ is given by
\begin{align*}
E[\xi_j]=\frac{1}{l(n)}+2 \sum_{h=1}^{l(n)} \frac{1}{l(n)}\left\{1-\frac{h}{l(n)}\right\}=1.
\end{align*}
For the variance, direct calculations yield
\begin{align*}
Var[\xi_j]=\left[\frac{1}{l(n)^2}+2 \sum_{h=1}^{l(n)}\frac{\{l(n)-h\}^2}{l(n)^4} \right]Var[w_.]=\left\{\frac{2}{3l(n)}+\frac{1}{3l(n)^3}\right\} Var[w_.]=1
\end{align*}
for all $j \in \mathbb{Z}.$ For any $j \in \mathbb{Z}$ and $|h|<2l(n)-1,$ the covariance function $\Cov(\xi_j,\xi_{j+h})$ can be described by a parabola centered at zero and opening downward \citep[for details, see ][Section 6.2]{buehlmann93}. The resulting sequence $(\xi_j)_{j\in \mathbb{Z}}$ satisfies \bf A1\it, \bf A2, \it and \bf A3\it. Figure \ref{fig:taperedmultipliers2} provides an illustration of the kernel function $\kappa_2$ as well as simulated trajectories of Rademacher-type tapered block multiplier random variables. In this setting, $(\xi_j)_{j\in \mathbb{Z}}$ satisfies \bf A1\it, \bf A2, \it and \bf A3b\it. Notice the smoothing which is driven by the choice of kernel function and the block length $l(n).$ This effect can be further explored using more sophisticated kernel functions, e.g., with bell-shape; this is left for further research.
\end{example}

\subsection{Finite sample behavior}
\label{subsection:NI_fs}

In the present section we investigate and compare the finite sample properties of the (moving) block bootstrap and the tapered block multiplier technique by means of a simulation study. For that purpose we use both resampling techniques in order to estimate the (co)variances of the empirical copula process.
This study complements the one of \citet{buecherdette09} and \citet{buecher11} on bootstrap approximations for the empirical copula process in the i.i.d setting.

To be more precise, in Tables \ref{tab:meanmcresults100}, \ref{tab:meanmcresults200} and \ref{tab:meanmcresultsARGARCH} we demonstrate simulation results (for $n=100$ and $n=200$) on the estimation of the theoretical covariance $\Cov(\GG_C(\mathbf{u}),\GG_C(\mathbf{v}))$ calculated at the points $\mathbf{u}=\mathbf{v} \in \{(1/3,1/3),(2/3,1/3),(1/3,2/3),(2/3,2/3)\}$ in basically six different models:
serial dependence features arise from (bivariate) i.i.d., AR($1$) or GARCH($1,1$) time series models, while the copulas linking the marginals of the innovations are taken from the Gumbel--Hougaard or the Clayton family. The technical details of the execution of the study are given below, we start with the discussion of the results.

The results based on independent and identically distributed samples indicate that the tapered block multiplier outperforms the block bootstrap in mean and MSE of estimation. However, applying the resampling methods of the present paper in an i.i.d. setting comes at the price of an slight increased mean squared error in comparison to the multiplier or (block) bootstrap with block length $l=1$ as investigated in \cite{buecherdette09}. Hence, we suggest to test serial independence of continuous multivariate time-series as introduced by \citet{kojadinovicyan11} to investigate which method is appropriate. In the case of serially dependent observations, the results indicate that the tapered block multiplier yields more precise results in mean and mean squared error than the block bootstrap (which tends to overestimate) for the considered choices of the temporal dependence structure, the kernel function and the copula. This finding coincides with the results in \cite{buecherdette09}.
Regarding the choice of the kernel function, mean results for $\kappa_1$ and $\kappa_2$ are similar, whereas $\kappa_2$ yields slightly better results in mean squared error. Additional MC simulations are given in \citet{ruppert11}: if the multiplier or bootstrap methods for independent observations are incorrectly applied to dependent observations, i.e., $l_B=l_M=1,$ then their results do not reflect the changed structure adequately. Results based on Normal, Gamma, and Rademacher-type sequences $(w_j)_{j \in \mathbb{Z}}$ indicate that different distributions used to simulate the multiplier random variables lead to similar results. To ease comparison of the next section with the work of \citet{remillardscaillet09}, we decided to state the results for normal multiplier random variables here.

\begin{table}[htbp!]
\caption{\small{\textbf{Mean and MSE $(\times 10^4)$ Monte Carlo results}. I.i.d. and AR($1$) settings, sample size $n=100$ and $1,000$ Monte Carlo replications.  For each replication, we perform $S=2,000$ tapered block multiplier $(M_i)$ repetitions with Normal multiplier random variables, kernel function $\kappa_i,$ $i=1,2,$ block length $l_M=3,$ and block bootstrap $(B)$ repetitions with block length $l_B=5.$}}
\label{tab:meanmcresults100}
\begin{center}
\small
\begin{tabular}{llrrrrrrrr} \hline
\\
$(u_1,u_2)$ &  & \multicolumn{2}{l}{$(1/3,1/3)$}  & \multicolumn{2}{l}{$(1/3,2/3)$} & \multicolumn{2}{l}{$(2/3,1/3)$} & \multicolumn{2}{l}{$(2/3,2/3)$}   \\
& & Mean & MSE & Mean & MSE & Mean & MSE & Mean & MSE \\
\hline \hline
\\
\multicolumn{10}{l}{\bf i.i.d. setting \rm} \\
Clayton & True  & 0.0486 & & 0.0338 & & 0.0338 & & 0.0508 &  \\
$(\theta=1)$ &  Approx. & 0.0487 & & 0.0338 & & 0.0338 & & 0.0508 & \\
& $M_2$  & 0.0496 & 1.6323 & 0.0344 & 1.2449 & 0.0345 & 1.3456 & 0.0528 & 1.6220 \\
& $M_1$  & 0.0494 & 1.7949 & 0.0342 & 1.2934 & 0.0343 & 1.4128 & 0.0524 & 1.7910 \\
& $B$ & 0.0599 & 2.8286 & 0.0432 & 2.4103 & 0.0429 & 2.1375 & 0.0643 & 3.1538 \\ \\

Clayton & True  & 0.0254 & & 0.0042 & & 0.0042 & & 0.0389 &  \\
$(\theta=4)$ &  Approx. & 0.0255 & & 0.0042 & & 0.0042 & & 0.0390 & \\
& $M_2$  & 0.0259 & 0.9785 & 0.0051 & 0.3715 & 0.0048 & 0.3662 & 0.0407 & 1.6324 \\
& $M_1$  & 0.0257 & 1.0104 & 0.0050 & 0.3656 & 0.0048 & 0.3673 & 0.0404 & 1.7048 \\
& $B$ & 0.0383 & 2.5207 & 0.0100 & 0.7222 & 0.0097 & 0.6662 & 0.0533 & 3.7255 \\ \\

Gumbel & True  & 0.0493 & & 0.0336 & & 0.0336 & & 0.0484 &  \\
$(\theta=1.5)$ &  Approx. & 0.0493 & & 0.0335 & & 0.0335 & & 0.0485 & \\
& $M_2$  & 0.0514 & 1.3914 & 0.0346 & 1.2657 & 0.0340 & 1.2583 & 0.0497 & 1.4946 \\
& $M_1$  & 0.0510 & 1.5530 & 0.0344 & 1.3390 & 0.0338 & 1.3275 & 0.0495 & 1.6948 \\
& $B$ & 0.0616 & 2.8334 & 0.0429 & 2.1366 & 0.0432 & 2.2273 & 0.0620 & 3.2134 \\ \\

Gumbel & True  & 0.0336 & & 0.0058 & & 0.0058 & & 0.0293 &  \\
$(\theta=3)$ &  Approx. & 0.0335 & & 0.0058 & & 0.0058 & & 0.0294 & \\
& $M_2$  & 0.0355 & 1.1359 & 0.0064 & 0.4427 & 0.0063 & 0.3851 & 0.0307 & 0.9819 \\
& $M_1$  & 0.0353 & 1.1991 & 0.0063 & 0.4409 & 0.0062 & 0.3832 & 0.0306 & 1.0261 \\
& $B$ & 0.0470 & 2.7885 & 0.0120 & 0.8396 & 0.0122 & 0.8959 & 0.0437 & 2.9355 \\ \\

\multicolumn{10}{l}{\bf AR($1$) setting with $\beta=0.5$ \rm} \\

Clayton &  Approx. & 0.0599 & & 0.0408 & & 0.0409 & & 0.0629 & \\
$(\theta=1)$ & $M_2$  & 0.0602 & 3.1797 & 0.0394 & 2.4919 & 0.0398 & 2.5903 & 0.0625 & 2.6297 \\
& $M_1$  & 0.0598 & 3.5305 & 0.0391 & 2.6090 & 0.0396 & 2.7410 & 0.0620 & 2.9826 \\
& $B$ & 0.0699 & 3.4783 & 0.0496 & 2.9893 & 0.0492 & 2.8473 & 0.0761 & 4.0660 \\ \\

Clayton &  Approx. & 0.0329 & & 0.0064 & & 0.0064 & & 0.0432 & \\
$(\theta=4)$ & $M_2$  & 0.0347 & 2.0017 & 0.0071 & 0.6040 & 0.0072 & 0.6666 & 0.0460 & 2.8656 \\
& $M_1$  & 0.0344 & 2.0672 & 0.0071 & 0.5869 & 0.0072 & 0.6583 & 0.0458 & 3.0571 \\
& $B$ & 0.0472 & 3.5289 & 0.0132 & 1.0990 & 0.0128 & 0.9658 & 0.0585 & 4.2965 \\ \\

Gumbel &  Approx. & 0.0617 & & 0.0408 & & 0.0409 & & 0.0605 & \\
$(\theta=1.5)$ & $M_2$  & 0.0631 & 3.0587 & 0.0406 & 2.7512 & 0.0398 & 2.4187 & 0.0600 & 3.0433 \\
& $M_1$  & 0.0626 & 3.4252 & 0.0402 & 2.8739 & 0.0395 & 2.4920 & 0.0594 & 3.3824 \\
& $B$ & 0.0735 & 3.7410 & 0.0501 & 3.2025 & 0.0502 & 3.1150 & 0.0735 & 4.0521 \\ \\

Gumbel &  Approx. & 0.0385 & & 0.0061 & & 0.0061 & & 0.0345 & \\
$(\theta=3)$ & $M_2$  & 0.0425 & 2.5441 & 0.0069 & 0.5841 & 0.0070 & 0.5796 & 0.0375 & 2.1518 \\
& $M_1$  & 0.0422 & 2.7239 & 0.0069 & 0.5706 & 0.0070 & 0.5836 & 0.0373 & 2.2381 \\
& $B$ & 0.0535 & 3.8803 & 0.0127 & 1.0720 & 0.0125 & 1.0526 & 0.0507 & 4.1894 \\ \\

\hline
\end{tabular}
\normalsize
\end{center}
\end{table}
\begin{table}[htbp!]
\caption{\small{\textbf{Mean and MSE $(\times 10^4)$ Monte Carlo results}. I.i.d. and AR($1$) settings, sample size $n=200$ and $1,000$ Monte Carlo replications.  For each replication, we perform $S=2,000$ tapered block multiplier $(M_i)$ repetitions with Normal multiplier random variables, kernel function $\kappa_i,$ $i=1,2,$ block length $l_M=4,$ and block bootstrap $(B)$ repetitions with block length $l_B=7.$}}
\label{tab:meanmcresults200}
\begin{center}
\small
\begin{tabular}{llrrrrrrrr} \hline
\\
$(u_1,u_2)$ &  & \multicolumn{2}{l}{$(1/3,1/3)$}  & \multicolumn{2}{l}{$(1/3,2/3)$} & \multicolumn{2}{l}{$(2/3,1/3)$} & \multicolumn{2}{l}{$(2/3,2/3)$}   \\
& & Mean & MSE & Mean & MSE & Mean & MSE & Mean & MSE \\
\hline \hline
\\
\multicolumn{10}{l}{\bf i.i.d. setting \rm} \\
Clayton & True  & 0.0486 & & 0.0338 & & 0.0338 & & 0.0508 &  \\
$(\theta=1)$ &  Approx. & 0.0487 & & 0.0338 & & 0.0338 & & 0.0508 & \\
& $M_2$  & 0.0494 & 1.2579 & 0.0346 & 0.8432 & 0.0343 & 0.8423 & 0.0522 & 1.0987 \\
& $M_1$  & 0.0490 & 1.3749 & 0.0345 & 0.8880 & 0.0341 & 0.8719 & 0.0519 & 1.2074 \\
& $B$ & 0.0562 & 1.7155 & 0.0402 & 1.2154 & 0.0395 & 1.2279 & 0.0593 & 1.8137 \\ \\

Clayton & True  & 0.0254 & & 0.0042 & & 0.0042 & & 0.0389 &  \\
$(\theta=4)$ &  Approx. & 0.0255 & & 0.0042 & & 0.0042 & & 0.0390 & \\
& $M_2$  & 0.0261 & 0.5811 & 0.0044 & 0.1889 & 0.0048 & 0.2027 & 0.0390 & 0.9932 \\
& $M_1$  & 0.0260 & 0.6024 & 0.0044 & 0.1876 & 0.0048 & 0.2008 & 0.0388 & 1.0702 \\
& $B$ & 0.0347 & 1.6381 & 0.0076 & 0.3197 & 0.0073 & 0.2908 & 0.0489 & 2.1425 \\ \\

Gumbel & True  & 0.0493 & & 0.0336 & & 0.0336 & & 0.0484 &  \\
$(\theta=1.5)$ &  Approx. & 0.0493 & & 0.0335 & & 0.0335 & & 0.0485 & \\
& $M_2$  & 0.0512 & 1.1513 & 0.0347 & 0.8157 & 0.0347 & 0.8113 & 0.0504 & 1.1697 \\
& $M_1$  & 0.0511 & 1.2909 & 0.0345 & 0.8653 & 0.0346 & 0.8692 & 0.0503 & 1.2811 \\
& $B$ & 0.0577 & 1.7178 & 0.0402 & 1.2233 & 0.0403 & 1.1701 & 0.0574 & 1.8245 \\ \\

Gumbel & True  & 0.0336 & & 0.0058 & & 0.0058 & & 0.0293 &  \\
$(\theta=3)$ &  Approx. & 0.0335 & & 0.0058 & & 0.0058 & & 0.0294 & \\
& $M_2$  & 0.0361 & 0.8340 & 0.0067 & 0.2577 & 0.0063 & 0.2505 & 0.0320 & 0.8678 \\
& $M_1$  & 0.0359 & 0.9107 & 0.0067 & 0.2576 & 0.0062 & 0.2444 & 0.0318 & 0.9078 \\
& $B$ & 0.0435 & 1.7448 & 0.0095 & 0.3917 & 0.0094 & 0.3804 & 0.0388 & 1.5320 \\ \\

\multicolumn{10}{l}{\bf AR($1$) setting with $\beta=0.5$ \rm} \\

Clayton &  Approx. & 0.0599 & & 0.0408 & & 0.0409 & & 0.0629 & \\
$(\theta=1)$ & $M_2$  & 0.0615 & 2.3213 & 0.0413 & 1.8999 & 0.0414 & 1.9235 & 0.0646 & 2.2534 \\
& $M_1$  & 0.0608 & 2.5553 & 0.0408 & 1.9594 & 0.0410 & 2.0094 & 0.0638 & 2.4889 \\
& $B$ & 0.0676 & 2.6206 & 0.0468 & 1.9802 & 0.0472 & 2.0278 & 0.0715 & 2.5494 \\ \\

Clayton &  Approx. & 0.0329 & & 0.0064 & & 0.0064 & & 0.0432 & \\
$(\theta=4)$ & $M_2$  & 0.0344 & 1.2685 & 0.0074 & 0.3890 & 0.0073 & 0.4108 & 0.0460 & 1.7936 \\
& $M_1$  & 0.0341 & 1.3034 & 0.0073 & 0.3824 & 0.0072 & 0.4053 & 0.0455 & 1.8566 \\
& $B$ & 0.0432 & 2.1693 & 0.0105 & 0.5427 & 0.0101 & 0.4901 & 0.0546 & 2.8808 \\ \\

Gumbel &  Approx. & 0.0617 & & 0.0408 & & 0.0409 & & 0.0605 & \\
$(\theta=1.5)$ & $M_2$  & 0.0649 & 2.2833 & 0.0432 & 2.1528 & 0.0433 & 1.9800 & 0.0646 & 3.0536 \\
& $M_1$  & 0.0639 & 2.4663 & 0.0426 & 2.1943 & 0.0428 & 2.0367 & 0.0638 & 3.2590 \\
& $B$ & 0.0703 & 2.5133 & 0.0468 & 1.8163 & 0.0467 & 1.7087 & 0.0693 & 2.5655 \\ \\

Gumbel &  Approx. & 0.0385 & & 0.0061 & & 0.0061 & & 0.0345 & \\
$(\theta=3)$ & $M_2$  & 0.0430 & 1.9156 & 0.0079 & 0.4884 & 0.0074 & 0.3950 & 0.0405 & 2.5846 \\
& $M_1$  & 0.0426 & 1.9784 & 0.0078 & 0.4806 & 0.0073 & 0.3858 & 0.0400 & 2.6083 \\
& $B$ & 0.0492 & 2.2152 & 0.0102 & 0.5375 & 0.0099 & 0.4832 & 0.0455 & 2.1534 \\ \\

\hline
\end{tabular}
\normalsize
\end{center}
\end{table}

\begin{table}[htbp!]
\caption{\small{\textbf{Mean and MSE $(\times 10^4)$ Monte Carlo results}. AR($1$) and GARCH($1,1$) settings, sample size $n=200$ and $1,000$ Monte Carlo replications. For each replication, we perform $S=2,000$ tapered block multiplier $(M_i)$ repetitions with Normal multiplier random variables, kernel function $\kappa_i,$ $i=1,2,$ block length $l_M=4,$ and block bootstrap $(B)$ repetitions with block length $l_B=7.$}}
\label{tab:meanmcresultsARGARCH}
\begin{center}
\small
\begin{tabular}{llrrrrrrrr} \hline
\\
$(u_1,u_2)$ &  & \multicolumn{2}{l}{$(1/3,1/3)$}  & \multicolumn{2}{l}{$(1/3,2/3)$} & \multicolumn{2}{l}{$(2/3,1/3)$} & \multicolumn{2}{l}{$(2/3,2/3)$}   \\
& & Mean & MSE & Mean & MSE & Mean & MSE & Mean & MSE \\
\hline \hline
\\

\multicolumn{10}{l}{\bf AR($1$) setting with $\beta=0.25$ \rm} \\
Clayton &  Approx. & 0.0506 & & 0.0350 & & 0.0350 & & 0.0530 & \\
$(\theta=1)$ & $M_2$  & 0.0521 & 1.5319 & 0.0360 & 1.0527 & 0.0361 & 1.1258 & 0.0545 & 1.4689 \\
& $M_1$  & 0.0517 & 1.6574 & 0.0357 & 1.1003 & 0.0358 & 1.1662 & 0.0541 & 1.6083 \\
& $B$ & 0.0591 & 2.0394 & 0.0415 & 1.4084 & 0.0417 & 1.5256 & 0.0624 & 2.1027 \\ \\

Clayton &  Approx. & 0.0272 & & 0.0048 & & 0.0048 & & 0.0395 & \\
$(\theta=4)$ & $M_2$  & 0.0285 & 0.8366 & 0.0052 & 0.2290 & 0.0052 & 0.2208 & 0.0413 & 1.3315 \\
& $M_1$  & 0.0283 & 0.8792 & 0.0051 & 0.2273 & 0.0052 & 0.2200 & 0.0410 & 1.3858 \\
& $B$ & 0.0375 & 1.9937 & 0.0084 & 0.3726 & 0.0085 & 0.3646 & 0.0495 & 2.1996 \\ \\

Gumbel &  Approx. & 0.0518 & & 0.0349 & & 0.0348 & & 0.0507 & \\
$(\theta=1.5)$ & $M_2$  & 0.0549 & 1.4150 & 0.0373 & 1.1973 & 0.0377 & 1.2459 & 0.0550 & 2.0067 \\
& $M_1$  & 0.0545 & 1.5154 & 0.0370 & 1.2339 & 0.0373 & 1.2871 & 0.0545 & 2.0766 \\
& $B$ & 0.0608 & 2.0500 & 0.0415 & 1.3234 & 0.0418 & 1.4093 & 0.0608 & 2.3118 \\ \\

Gumbel &  Approx. & 0.0346 & & 0.0058 & & 0.0058 & & 0.0304 & \\
$(\theta=3)$ & $M_2$  & 0.0386 & 1.1686 & 0.0070 & 0.3079 & 0.0068 & 0.2979 & 0.0350 & 1.5135 \\
& $M_1$  & 0.0384 & 1.2224 & 0.0070 & 0.3052 & 0.0067 & 0.2892 & 0.0347 & 1.5157 \\
& $B$ & 0.0447 & 1.8434 & 0.0095 & 0.3873 & 0.0097 & 0.4303 & 0.0409 & 1.8797 \\ \\

\multicolumn{10}{l}{\bf GARCH($1,1$) setting \rm} \\
Clayton &  Approx. & 0.0479 & & 0.0340 & & 0.0340 & & 0.0516 & \\
$(\theta=1)$ & $M_2$  & 0.0491 & 1.1144 & 0.0347 & 0.8485 & 0.0343 & 0.8279 & 0.0520 & 1.0958 \\
& $M_1$  & 0.0486 & 1.3579 & 0.0339 & 0.9021 & 0.0338 & 0.8100 & 0.0515 & 1.2013 \\
& $B$ & 0.0567 & 2.0156 & 0.0403 & 1.1765 & 0.0403 & 1.2556 & 0.0600 & 1.8542 \\ \\

Clayton &  Approx. & 0.0252 & & 0.0055 & & 0.0056 & & 0.0403 & \\
$(\theta=4)$ & $M_2$  & 0.0259 & 0.5054 & 0.0051 & 0.1979 & 0.0053 & 0.2301 & 0.0399 & 1.0431 \\
& $M_1$  & 0.0258 & 0.6429 & 0.0052 & 0.2199 & 0.0051 & 0.2217 & 0.0390 & 1.0959 \\
& $B$ & 0.0345 & 1.5252 & 0.0081 & 0.2764 & 0.0081 & 0.2921 & 0.0484 & 1.8359 \\ \\

Gumbel & Approx. & 0.0500 & & 0.0339 & & 0.0339 & & 0.0482 & \\
$(\theta=1.5)$ & $M_2$  & 0.0516 & 1.0480 & 0.0356 & 0.8486 & 0.0354 & 0.8175 & 0.0511 & 1.2451 \\
& $M_1$  & 0.0516 & 1.2235 & 0.0351 & 0.9582 & 0.0352 & 0.8848 & 0.0503 & 1.2774 \\
& $B$ & 0.0575 & 1.5928 & 0.0402 & 1.2395 & 0.0403 & 1.2346 & 0.0574 & 1.9198 \\ \\

Gumbel &  Approx. & 0.0341 & & 0.0074 & & 0.0074 & & 0.0291 & \\
$(\theta=3)$ & $M_2$  & 0.0362 & 0.9284 & 0.0073 & 0.2941 & 0.0072 & 0.2587 & 0.0321 & 0.9133 \\
& $M_1$  & 0.0366 & 0.9782 & 0.0073 & 0.2786 & 0.0071 & 0.2888 & 0.0320 & 0.9819 \\
& $B$ & 0.0435 & 1.6811 & 0.0103 & 0.3607 & 0.0101 & 0.3390 & 0.0390 & 1.6878 \\ \\

\hline
\end{tabular}
\normalsize
\end{center}
\end{table}

Finally, for the sake of completeness we state the technical details underlying the simulation. First of all, the details regarding the models are as follows.
\begin{itemize}
\item
The bivariate Clayton and Gumbel--Hougaard copulas are given by
\begin{align*}
C_\theta^{Cl}(u_1,u_2)&=\left(u_1^{-\theta}+u_2^{-\theta}-1\right)^{-\frac{1}{\theta}}, \ \theta > 0, \\
C_\theta^{Gu}(u_1,u_2)&=\exp\left(-\left[\left\{-\ln\left(u_1\right)\right\}^\theta+\left\{-\ln\left(u_2\right)\right\}^\theta\right]^{\frac{1}{\theta}}\right), \ \theta \geq 1,
\end{align*}
respectively, and we chose the parameters in such a way that Kendall's $\tau$ is either $1/3$ or $2/3$, i.e., $\theta \in \{1, 4\}$ for the Clayton and $\theta \in \{1.5, 3\}$ for the Gumbel--Hougaard copula.

\item
The AR(1) process we consider is the stationary solution of the AR(1) equation
\[
	\mathbf{X}_j=\beta \mathbf{X}_{j-1}+\boldsymbol{\varepsilon}_j,
\]
where the innovations are supposed to have standard normal marginals linked by one of the aforementioned copulas and where the coefficient of the lagged variable is either $\beta=0.25$ or $\beta=0.5$. This stationary solution can be written as $\mathbf{X}_j=\sum_{k=0}^\infty \beta^{k} \boldsymbol{\eps}_{j-k}$.
 We simulate an (approximate) sample of length $n$ of this model as follows:  for some reasonably large negative number $M$, e.g., $M=-100$, let $\mathbf{U}_j=(U_{j,1},U_{j,2})$, $j=M,\dots,n$ be a sample of independent realizations of one of the aforementioned copulas. Set $\boldsymbol{\varepsilon}_j=(\Phi^{-1}(U_{j,1}), \Phi^{-1}(U_{j,2}))$, with $\Phi$ being the standard normal cdf, and recursively define $\mathbf{X}_M=\boldsymbol{\varepsilon}_M$ and
\begin{align}
\label{equation:AR_MC_3}
\mathbf{X}_j=\beta \mathbf{X}_{j-1}+\boldsymbol{\varepsilon}_j \ \text{for all} \ j=M+1, \ldots, n.
\end{align}
The last $n$ observations form the sample $\mathbf{X}_1,\dots,\mathbf{X}_n$.

\item
The GARCH(1,1) sample is simulated as following: define $\boldsymbol{\varepsilon}_j=(\eps_{j,1},\eps_{j,2})$ as in the AR(1)-example and recursively define
$\sigma_{M,i}=\sqrt{{\omega_i}/({1-\alpha_i-\beta_i})} \ \text{for} \ i=1,2$
and
\begin{align}
\label{equation:GARCH1}
X_{j,i}=\sigma_{j,i} \eps_{j,i}, \text{ with } \sigma_{j,i}^2=\omega_i + \beta_i \sigma_{j-1,i}^2+\alpha_i \eps_{j-1,i}^2
\end{align}
for $i=1,2$ and $j=M+1,\dots, n$, where $\omega_1=0.012, \omega_2=0.037, \beta_1=0.919, \beta_2=0.868, \alpha_1=0.072$ and $\alpha_2=0.115$. Again, the last $n$ observations form the sample $\mathbf{X}_1,\dots,\mathbf{X}_n$.
The considered coefficients are estimates derived in \citet{jondeauetal07} to model volatility of S\&P $500$ and DAX daily (log-)returns in an empirical application which shows the practical relevance of this specific parameter choice.
\end{itemize}

The theoretical covariance of $\GG_C$ is easily calculated in the i.i.d. setting, but hardly derivable in the serially dependent case. Note that even a closed form expression for the copula of $\mathbf{X}_j$ is unaccessible. For that reason, we approximate the theoretical covariance by means of empirical covariances calculated from simulated samples of the process $\sqrt{n} (\widehat{C}_n-\widehat{C}_N)$, where $n=n(N)\to \infty$ with $n=o(N)$. We chose $N=10^6$ and $n=1000$, and made $10^6$ replications on basis of which we calculated the empirical covariances. The results in Tables  \ref{tab:meanmcresults100} - \ref{tab:meanmcresultsARGARCH} (the `True' vs. `Approx.' lines) for the i.i.d. setting show that this approximation works sufficiently well, whence we can use it as a benchmark in the two serial dependent settings.

Regarding the tapered block multiplier bootstrap  we decided to use both kernel functions $\kappa_1$ and $\kappa_2$ from Examples  \ref{example:uniform_blocks} and \ref{example:triangular_blocks}. For the (moving) block bootstrap we choose the block length as $l_B(100)=5$ and $l_B(200)=7;$ this choice corresponds to $l_B(n)=\lfloor 1.25n^{1/3} \rfloor$ which satisfies the assumptions of the asymptotic theory. For a detailed discussion on the block length of the block bootstrap, we refer to \citet{kuensch89} as well as \citet{buehlmannkuensch99}.
The tapered block multiplier technique is assessed based on a sequence $(w_j)_{j \in \mathbb{Z}}$ of normal random variables as introduced in Examples \ref{example:uniform_blocks} and \ref{example:triangular_blocks}.
The block length is set to $l_M(n)=\lfloor 1.1n^{1/4} \rfloor,$ hence $l_M(100)=3$ and $l_M(200)=4,$ meaning that both methods yield $2l_M$-dependent blocks. For all methods $S=2,000$ bootstrap repetitions are performed and the target covariance is estimated by the sample covariance over the 2,000 repetitions.
This procedure is repeated $1,000$ times and we report mean and mean squared error (MSE) for each method.

\section{Testing for a constant copula}
\label{section:tests}

The present section presents two nonparametric tests for a constant copula in the case of serially dependent processes. We begin with a test where a change point candidate is given, and proceed with a more general test without specifying a time point where a break in the copula structure occurs. Both tests are
 consistent against general alternatives and their finite sample performance is investigated by means of a simulation study.

\subsection{Specified change point candidate}
\index{test for a constant copula! specified change point candidate}

The specification of a change point candidate can for instance have an economic motivation: \citet{patton02} investigates a change in parameters of the dependence structure between various exchange rates following the introduction of the euro on the $1$st of January $1999.$ Focusing on stock returns, multivariate association between major S\&P global sector indices before and after the bankruptcy of Lehman Brothers Inc. on $15$th of September $2008$ is assessed in \citet{gaisseretal10} and \citet{ruppert11}. Whereas these references investigate change points in functionals of the copula, the copula itself is in the focus of this study. This approach permits to analyze changes in the structure of association even if a functional thereof, such as a measure of multivariate association, is invariant.

Suppose we observe a sample $\mathbf{X}_1, \ldots, \mathbf{X}_n$ of a process $(\mathbf{X}_j)_{j \in \mathbb{Z}}.$ We derive a test for constancy of the copula in the case of a specified change point candidate indexed by $\lfloor \lambda n \rfloor$ for $\lambda \in (0,1)$, i.e., a test for the hypothesis
\begin{align*}
& H_0: \mathbf{U_j} \sim C_1 \ \text{for all} \ j=1, \ldots, n, \\
& H_1: \mathbf{U_j} \sim \begin{cases}   C_1 \ \text{for all} \ j=1, \ldots, \lfloor \lambda n \rfloor, \\ C_2 \ \text{for all} \ j=\lfloor \lambda n \rfloor +1, \ldots, n, \end{cases}
\end{align*}
where $C_1$ and $C_2$ are assumed to be different in at least one point $\mathbf u \in [0,1]^d$ (and hence, for continuity reasons, also in a neighbourhood of this point).  The proposed test statistic will be based on a splitting of the sample into two subsamples: $\mathbf{X}_1, \ldots, \mathbf{X}_{\lfloor \lambda n \rfloor}$ and $\mathbf{X}_{\lfloor \lambda n \rfloor+1},\ldots,\mathbf{X}_n.$ A significant discrepancy between estimates of the copula in the two subsamples suggests to reject the null hypothesis.
Assuming constant marginal distributions in each subsample, let
\begin{align*}
	& \widehat{C}_{1,\dots,\lfloor \lambda n \rfloor}(\mathbf{u}) = \frac{1}{\lfloor \lambda n \rfloor}\sum_{j=1}^{\lfloor \lambda n \rfloor}  \mathbf{1}_{\{\widehat{\mathbf{U}}_{j} \leq \mathbf{u}\}}, &
	& \widehat{C}_{\lfloor \lambda n \rfloor + 1, \dots, n}(\mathbf{u}) = \frac{1}{n-\lfloor \lambda n \rfloor}\sum_{j=\lfloor \lambda n \rfloor+1}^n  \mathbf{1}_{\{\widehat{\mathbf{V}}_{j} \leq \mathbf{u}\}} &
\end{align*}
denote the corresponding empirical copulas of the two subsamples, where $\widehat{\mathbf{U}}_{1}, \dots, \widehat{\mathbf{U}}_{\lfloor \lambda n \rfloor}$ and $\widehat{\mathbf{V}}_{\lfloor \lambda n \rfloor+1}, \dots, \widehat{\mathbf{V}}_{n}$
denote the pseudo-observations calculated from the first and second subsample, respectively.
We consider the test statistic defined by
\begin{align}
\label{equation:stat_specbreak}
 	T_n(\lambda)= \int_{[0,1]^d}\left[\sqrt{\frac{\lfloor \lambda n \rfloor(n-\lfloor \lambda n \rfloor)}{n}}\left\{\widehat{C}_{1,\dots, \lfloor \lambda n \rfloor}(\mathbf{u})-\widehat{C}_{\lfloor \lambda n \rfloor + 1,\dots, n}(\mathbf{u})\right\}\right]^2 d\mathbf{u}
\end{align}
which can be calculated explicitly; for details we refer to \citet{remillardscaillet09}. These authors introduce a test for equality between two copulas which is applicable in the case of serial independence. Weak convergence of $T_n(\lambda)$ under strong mixing follows from the following result.

\begin{theorem}
\label{theorem:asymptoticsS}
Consider observations $\mathbf{X}_{1},\ldots,\mathbf{X}_{n},$ drawn from a process $(\mathbf{X}_j)_{j \in \mathbb{Z}}$ satisfying the strong mixing condition $\alpha_\mathbf{X}(r)=\mathcal{O}(r^{-a})$ for some $a>1.$ Further assume a specified change point candidate indexed by $\lfloor \lambda n \rfloor$ for $\lambda \in (0,1)$ such that $\mathbf{U_j} \sim C_1,$ $X_{j,i} \sim F_{1,i}$ for all $j=1, \ldots, \lfloor \lambda n \rfloor,$ $i=1,\ldots,d$ and $\mathbf{U_j} \sim C_2,$ $X_{j,i} \sim F_{2,i}$ for all $j=\lfloor \lambda n \rfloor+1, \ldots,n,$ $i=1,\ldots,d.$ Suppose that $C_1$ and $C_2$ satisfy Condition (\ref{equation:segers}). Under the null hypothesis $C_1=C_2$, in the metric space $\ell^\infty([0,1]^d)$,
\[
	\sqrt{\frac{\lfloor \lambda n \rfloor(n-\lfloor \lambda n \rfloor)}{n}}\left\{\widehat{C}_{1, \dots,\lfloor \lambda n \rfloor} - \widehat{C}_{\lfloor \lambda n \rfloor+1, \dots, n}\right\} \weak  \sqrt{1-\lambda} \GG_{C_1} - \sqrt{\lambda} \GG_{C_2},
\]
where $\GG_{C_p}(\mathbf{u}) = \BB_{C_p} (\mathbf{u}) - \sum_{i=1}^d D_iC(\mathbf{u}) \BB_{C_p}(\mathbf{u}^{(i)})$ for $p=1,2$ and where $\BB_{C_1}$ and $\BB_{C_2}$ denote two independent tight centered Gaussian processes with covariances as specified in~\eqref{equation:covariancestructure}.
\end{theorem}

As a consequence of this Theorem, the continuous mapping Theorem yields
\begin{align*}
	T_n(\lambda) \weak T(\lambda)=\int_{[0,1]^d} \left\{ \sqrt{1-\lambda} \GG_{C_1}(\mathbf{u}) - \sqrt{\lambda} \GG_{C_2}(\mathbf{u}) \right\}^2 d\mathbf{u}
\end{align*}
for all $\lambda\in(0,1)$.
Note that if there exists a subset $\mathcal{I} \in [0,1]^d$ such that
\begin{align*}
\int_{\mathcal{I}} \left\{ \sqrt{1-\lambda} C_1(\mathbf{u}) - \sqrt{\lambda} C_2(\mathbf{u}) \right\}^2 d\mathbf{u}>0,
\end{align*}
then $T_n(\lambda) \rightarrow \infty$ in probability under $H_1.$
To estimate p-values of the test statistic, we use the tapered block multiplier technique developed in Section \ref{section:inference}. For that purpose, let $(\xi_{n,j})_{j \in \mathbb{Z}}$ be a sequence of multipliers satisfying \bf A1\rm, \bf A2\rm, \bf A3 \rm  and set
\begin{align*}
	\widehat{\BB}_{C,n,1}^M (\mathbf{u}) &= \frac{1}{\sqrt{\lfloor \lambda n\rfloor}}\sum_{j=1}^{\lfloor \lambda n\rfloor} \left( \frac{\xi_{n,j}}{\bar \xi_n^{(1)}} -1 \right) \mathbf{1}_{\{ \widehat{U}_j \le \mathbf{u} \}} \\
	\widehat{\BB}_{C,n,2}^M (\mathbf{u}) &= \frac{1}{\sqrt{n-\lfloor \lambda n\rfloor}}\sum_{j=\lfloor \lambda n\rfloor+1}^n \left( \frac{\xi_{n,j}}{\bar \xi_{n}^{(2)}} -1 \right) \mathbf{1}_{\{ \widehat{V}_j \le \mathbf{u} \}},
\end{align*}
where $\bar \xi_n^{(1)}$ and $\bar \xi_n^{(2)}$ denote the arithmetic mean of the $\xi_{j,n}$ in the corresponding samples.
Let $\widehat{D_iC_{n,1}}$ and $\widehat{D_iC_{n,2}}$ denote estimators for the corresponding partial derivatives and set, for $p=1,2$,
\[
	\widehat{\GG}_{C,n,p}^M(\mathbf{u}) = \widehat{\BB}_{C,n,p}^M(\mathbf{u}) - \sum_{i=1}^d \widehat{D_iC_{n,p}}(\mathbf{u}) \widehat{\BB}_{C,n,p}^M (\mathbf{u}^{(i)}).
\]
In the i.i.d.\ multiplier case this statistic is, up to some negligible constants, the same as the one used in \cite{remillardscaillet09} to test for equality between two copulas.

\begin{proposition}
 \label{prop:estimationS}
Consider observations $\mathbf{X}_{1},\ldots,\mathbf{X}_{\lfloor \lambda n \rfloor}$ and $\mathbf{X}_{\lfloor \lambda n \rfloor +1},\ldots,\mathbf{X}_{n}$ drawn from a process $(\mathbf{X}_j)_{j \in \mathbb{Z}}.$ Assume that the process satisfies the strong mixing assumptions of Theorem~\ref{theorem:taperedmultiplier}. Let $\lfloor \lambda n \rfloor$ for $\lambda \in (0,1)$ denote a specific change point candidate such that $\mathbf{U_j} \sim C_1,$ $X_{j,i} \sim F_{1,i}$ for all $j=1, \ldots, \lfloor \lambda n \rfloor,$ $i=1,\ldots,d$ and $\mathbf{U_j} \sim C_2,$ $X_{j,i} \sim F_{2,i}$ for all $j=\lfloor \lambda n \rfloor+1, \ldots,n,$ $i=1,\ldots,d.$ Suppose that $C_1$ and $C_2$ satisfy Condition (\ref{equation:segers}) and that $\widehat{D_iC_{n,1}}$ and $\widehat{D_iC_{n,2}}$ are estimators for the corresponding partial derivatives satisfying \bf C1\it and \bf C2\it. Let $\xi_{n,1}, \ldots, \xi_{n,n}$ denote samples of a tapered block multiplier process $(\xi_{n,j})_{j \in \mathbb{Z}}$ satisfying \bf A1\it, \bf A2\it, \bf A3 \it with block length $l(n) \rightarrow \infty,$ where $l(n)=\mathcal{O}(n^{1/2-\eps})$ for $0< \eps <1/2.$ Then
\begin{align}
\label{equation:tapmult_specbreak}
	\widehat{\HH}_{n,\lambda}^{M} := \sqrt{1-\lambda} \widehat{\GG}_{C,n,1}^{M} - \sqrt{\lambda} \widehat{\GG}_{C,n,2}^{M}
		\weakcondp{\xi} \sqrt{1-\lambda} \GG_{C_1} - \sqrt{\lambda} \GG_{C_2}
\end{align}
in $\ell^\infty([0,1]^d)$ both under the null hypothesis as well as under the alternative.
\end{proposition}

As a consequence, by the continuous mapping Theorem for the bootstrap, see \cite{kosorok08},
\begin{align*}
	\widehat{T}_n^M(\lambda) = \int_{[0,1]^d} (\widehat{\HH}_{n,\lambda}^M(\mathbf{u}))^2 d\mathbf{u} \weakcondp{\xi} T(\lambda).
\end{align*}
The integral involved in $\widehat{T}_n^{M}(\lambda)$ can be calculated explicitly \citep[see][Appendix B]{remillardscaillet06}. If we repeat the procedure $S$ times to obtain a sample $\widehat{T}_n^{M(1)}, \dots, \widehat{T}_n^{M(S)}$, then an approximate p-value for the test for $H_0$ is provided by
\begin{align}
\label{equation:estimate_p}
\frac{1}{S} \sum_{s=1}^S \mathbf{1}_{\left\{\widehat{T}_n^{M(s)}(\lambda)>T_n(\lambda)\right\}}.
\end{align}
Hence, p-values can be estimated by counting the number of cases in which the simulated test statistic based on the tapered block multiplier method exceeds the observed one.

\bigskip

\bf Finite sample properties. \rm Size and power of the test in finite samples are assessed in a simulation study. We consider bivariate samples of size $n=100$ or $n=200$ generated as in Section~\ref{subsection:NI_fs}, i.e., marginal i.i.d., AR(1) and GARCH(1,1) processes are either linked by a Clayton or a Gumbel--Hougaard copula. The change point after observation $\lfloor \lambda n \rfloor=n/2$ only affects the parameter within each family: the copula $C_1$ is parameterized such that Kendall's $\tau_1=0.2$, the copula $C_2$ such that Kendall's $\tau_2=0.2, \ldots, 0.9$. A set of $S=2,000$ normal tapered block multiplier processes is simulated, where the kernel function is chosen as $\kappa_2$ as suggested by the simulation results in Section \ref{subsection:NI_fs} and where $l_M(100)=3$ and $l_M(200)=4$ are chosen for the block length.
\begin{table}[t]
\caption{{\small{\bf Size and power of the test for a constant copula with a specified change point candidate}. Results are based on $1,000$ Monte Carlo replications, $n=100,$ $S=2,000$ tapered block multiplier repetitions, kernel function $\kappa_2,$ and asymptotic significance level $\alpha=5\%$.}}
\label{tab:powerkbn100}
\begin{center}
\small
\begin{tabular}{llrrrrrrrr} \hline
\\
$\tau_2$ & & $0.2$ & $0.3$ & $0.4$ & $0.5$ & $0.6$ & $0.7$ & $0.8$ & $0.9$ \\
\hline \hline
\vspace{-5pt} \\
\multicolumn{10}{l}{\bf i.i.d. setting \rm} \\
Clayton & $l=1$ & 0.036 & 0.110 & 0.295 & 0.612 & 0.881 & 0.983 & 1.000 & 1.000 \\
& $l=3$ & 0.050 & 0.114 & 0.315 & 0.578 & 0.877 & 0.976 & 1.000 & 1.000 \\
Gumbel & $l=1$ & 0.040 & 0.093 & 0.236 & 0.569 & 0.840 & 0.976 & 0.998 & 1.000 \\
& $l=3$ & 0.063 & 0.110 & 0.276 & 0.594 & 0.866 & 0.983 & 1.000 & 1.000 \\
\vspace{-5pt} \\
\multicolumn{10}{l}{\bf GARCH($1,1$) setting \rm} \\
Clayton & $l=1$ & 0.037 & 0.106 & 0.298 & 0.598 & 0.868 & 0.977 & 1.000 & 1.000 \\
& $l=3$ & 0.047 & 0.120 & 0.303 & 0.588 & 0.876 & 0.978 & 0.999 & 1.000 \\
Gumbel & $l=1$ & 0.043 & 0.089 & 0.246 & 0.573 & 0.827 & 0.978 & 0.999 & 1.000 \\
& $l=3$ & 0.065 & 0.124 & 0.285 & 0.569 & 0.847 & 0.980 & 1.000 & 1.000 \\
\vspace{-5pt} \\
\multicolumn{10}{l}{\bf AR($1$) setting with $\beta=0.25$ \rm} \\
Clayton & $l=1$ & 0.051 & 0.115 & 0.308 & 0.592 & 0.849 & 0.969 & 0.999 & 1.000 \\
& $l=3$ & 0.047 & 0.111 & 0.292 & 0.547 & 0.836 & 0.968 & 0.998 & 1.000 \\
Gumbel & $l=1$ & 0.053 & 0.109 & 0.257 & 0.550 & 0.836 & 0.975 & 0.998 & 1.000 \\
& $l=3$ & 0.066 & 0.105 & 0.254 & 0.568 & 0.818 & 0.964 & 1.000 & 1.000 \\
\vspace{-5pt} \\
\multicolumn{10}{l}{\bf AR($1$) setting with $\beta=0.5$ \rm} \\
Clayton & $l=1$ & 0.086 & 0.154 & 0.313 & 0.549 & 0.798 & 0.928 & 0.985 & 1.000 \\
& $l=3$ & 0.078 & 0.117 & 0.236 & 0.462 & 0.730 & 0.868 & 0.986 & 0.999 \\
Gumbel & $l=1$ & 0.100 & 0.172 & 0.285 & 0.541 & 0.816 & 0.956 & 0.998 & 1.000 \\
& $l=3$ & 0.077 & 0.109 & 0.218 & 0.482 & 0.722 & 0.907 & 0.994 & 0.999 \\
\\
\hline
\vspace{-20pt}
\end{tabular}
\normalsize
\end{center}
\end{table}

\begin{table}[t]
\caption{{\small{\bf Size and power of the test for a constant copula with a specified change point candidate}. Results are based on  $1,000$ Monte Carlo replications, $n=200,$ $S=2,000$ tapered block multiplier repetitions, kernel function $\kappa_2,$ and asymptotic significance level $\alpha=5\%$.}}
\label{tab:powerkbn200}
\begin{center}
\small
\begin{tabular}{llrrrrrrrr} \hline
\\
$\tau_2$ & & $0.2$ & $0.3$ & $0.4$ & $0.5$ & $0.6$ & $0.7$ & $0.8$ & $0.9$ \\
\hline \hline \\
\multicolumn{10}{l}{\bf i.i.d. setting \rm} \\
Clayton & $l=1$ & 0.047 & 0.172 & 0.524 & 0.908 & 0.993 & 1.000 & 1.000 & 1.000 \\
& $l=4$ & 0.063 & 0.164 & 0.552 & 0.905 & 0.989 & 1.000 & 1.000 & 1.000 \\
Gumbel & $l=1$ & 0.043 & 0.162 & 0.525 & 0.877 & 0.991 & 1.000 & 1.000 & 1.000 \\
& $l=4$ & 0.055 & 0.169 & 0.535 & 0.895 & 0.996 & 1.000 & 1.000 & 1.000 \\
\\
\multicolumn{10}{l}{\bf GARCH($1,1$) setting \rm} \\
Clayton & $l=1$ & 0.040 & 0.169 & 0.503 & 0.894 & 0.994 & 1.000 & 1.000 & 1.000 \\
& $l=4$ & 0.056 & 0.160 & 0.541 & 0.903 & 0.992 & 1.000 & 1.000 & 1.000 \\
Gumbel & $l=1$ & 0.046 & 0.154 & 0.498 & 0.873 & 0.992 & 1.000 & 1.000 & 1.000 \\
& $l=4$ & 0.057 & 0.174 & 0.496 & 0.899 & 0.994 & 1.000 & 1.000 & 1.000 \\
\\
\multicolumn{10}{l}{\bf AR($1$) setting with $\beta=0.25$ \rm} \\
Clayton & $l=1$ & 0.052 & 0.180 & 0.521 & 0.866 & 0.989 & 1.000 & 1.000 & 1.000 \\
& $l=4$ & 0.057 & 0.149 & 0.497 & 0.867 & 0.988 & 1.000 & 1.000 & 1.000 \\
Gumbel & $l=1$ & 0.047 & 0.180 & 0.515 & 0.872 & 0.989 & 1.000 & 1.000 & 1.000 \\
& $l=4$ & 0.050 & 0.136 & 0.490 & 0.855 & 0.992 & 1.000 & 1.000 & 1.000 \\
\\
\multicolumn{10}{l}{\bf AR($1$) setting with $\beta=0.5$ \rm} \\
Clayton & $l=1$ & 0.107 & 0.237 & 0.523 & 0.813 & 0.975 & 0.999 & 1.000 & 1.000 \\
& $l=4$ & 0.058 & 0.137 & 0.396 & 0.748 & 0.957 & 0.998 & 1.000 & 1.000 \\
Gumbel & $l=1$ & 0.122 & 0.227 & 0.499 & 0.825 & 0.979 & 1.000 & 1.000 & 1.000 \\
& $l=4$ & 0.059 & 0.123 & 0.401 & 0.770 & 0.958 & 0.998 & 1.000 & 1.000 \\
\\
\hline
\end{tabular}
\normalsize
\end{center}
\end{table}

The results of $1,000$ MC replications are shown in Tables \ref{tab:powerkbn100} and \ref{tab:powerkbn200} for $n=100$ and $n=200$, respectively. The test based on the tapered block multiplier technique  leads to a rejection quota under the null hypothesis which is close to the chosen theoretical asymptotic size of $5\%$ in all considered settings. Comparing the results for $n=100$ and $n=200,$ we observe that the approximation of the asymptotic size based on the tapered block multiplier improves in precision with increased sample size. The tapered block multiplier-based test also performs well under the alternative hypothesis and its power increases with the difference $\tau_2-\tau_1$ between the considered values for Kendall's $\tau.$ The power of the test under the alternative hypothesis is best in the case of no serial dependence as is shown in Table \ref{tab:powerkbn200}. If serial dependence is present in the sample then more observations are required to reach the power of the test in the case of serially independent observations. For comparison, we also show the results if the test assuming independent observations (i.e., the test based on the multiplier technique with block length $l=1$) is erroneously applied to the simulated dependent observations. The effects of different types of dependent observations differ largely in the finite sample simulations considered: GARCH($1,1$) processes do not show strong impact, whereas AR($1$) processes lead to considerable distortions, in particular regarding the size of the test. Results indicate that the test overrejects if temporal dependence is not taken into account; the observed size of the test in these cases can be more than twice the specified asymptotic size. For comparison, results for $n=200$ and kernel function $\kappa_1$ are shown in \citet{ruppert11}. The obtained results indicate that the uniform kernel function $\kappa_1$ leads to a more conservative testing procedure since the rejection quota is slightly higher, both under the null hypothesis as well as under the alternative. Due to the fact that the size of the test is approximated more accurately based on the kernel function $\kappa_2,$ its use is recommended.

\subsection{The general case: unspecified change point candidate}
\index{test for a constant copula! unspecified change point candidate}


The assumption of a change point candidate at specified location is relaxed in the following. Intuitively, testing with unspecified change point candidate(s) is less restrictive but a trade-off is to be made: the tests introduced in this section neither require conditions on the partial derivatives of the underlying copula(s) nor the specification of change point candidate(s), yet they are based on the assumption of strictly stationary univariate processes, i.e., $X_{j,i} \sim F_i$ for all $j\in \mathbb{Z}$ and $i=1,\ldots, d.$ The motivation for this test setting is that only for a subset of the change points documented in empirical studies, a priori hypothesis such as triggering economic events can be found \citep[see, e.g.,][]{diasembrechts09}. Even if a triggering event exists, its start (and end) often are subject to uncertainty: \citet{rodriguez06} studies changes in dependence structures of stock returns during periods of turmoil considering data framing the East Asian crisis in $1997$ as well as the Mexican devaluation in $1994,$ where no change point candidate is given a priori. These objects of investigation are well-suited for nonparametric methods which offer the important advantage that their results do not depend on model assumptions. For a general introduction to change point problems of this type, we refer to the monographs by \citet{csorgohorvath97} and, with particular emphasis on nonparametric methods, to \citet{brodskydarkhovsky93}.

Let $\mathbf{X}_1, \ldots, \mathbf{X}_n$ denote a sample of a process $(\mathbf{X}_j)_{j \in \mathbb{Z}}$ with strictly stationary univariate margins, i.e., $X_{j,i} \sim F_i$ for all $j\in \mathbb{Z}$ and $i=1,\ldots, d.$ We establish tests for the null hypothesis of a constant copula versus the alternative that there exist $P$ unspecified change points $\lambda_1<\ldots<\lambda_P \in [0,1],$ formally
\begin{align*}
H_0: \ & \mathbf{U_j} \sim C_1 \ \text{for all} \ j=1, \ldots, n, \\
H_1: \ & \text{there exist} \ 0=\lambda_0< \lambda_1 < \ldots < \lambda_{P}< \lambda_{P+1}=1 \ \text{such that} \ \mathbf{U_j} \sim C_p \\
& \text{for all} \ j=\lfloor \lambda_{p-1} n \rfloor +1, \ldots, \lfloor \lambda_p n \rfloor \ \text{and} \ p=1, \ldots, P+1,
\end{align*}
where, under the alternative hypothesis, $C_1, \ldots, C_{P+1}$ are assumed to be pairwise different in at least one point $\mathbf u \in [0,1]^d.$ Unlike in the previous section, we estimate the pseudo-observations $\widehat{\mathbf{U}}_1, \ldots, \widehat{\mathbf{U}}_n$ based on the whole sample $\mathbf{X}_1, \ldots, \mathbf{X}_n.$
The following test statistics are based on a comparison of the empirical distribution functions of the subsamples $\widehat{\mathbf{U}}_1, \ldots, \widehat{\mathbf{U}}_{\lfloor \zeta n \rfloor}$ and $\widehat{\mathbf{U}}_{\lfloor \zeta n \rfloor +1}, \ldots, \widehat{\mathbf{U}}_n$:
\begin{align*}
	\Sb_n(\zeta,\mathbf{u}) := \frac{\lfloor \zeta n \rfloor (n-\lfloor \zeta n \rfloor)}{n^{3/2}}
		\left\{ \frac{1}{\lfloor \zeta n \rfloor} \sum_{j=1}^{\lfloor \zeta n \rfloor}  \mathbf{1}_{\{\widehat{\mathbf{U}}_j \leq \mathbf{u}\}}  -
			\frac{1}{n-\lfloor \zeta n \rfloor} \sum_{j=\lfloor \zeta n \rfloor + 1}^n  \mathbf{1}_{\{\widehat{\mathbf{U}}_j \leq \mathbf{u}\}} \right\}.
\end{align*}
Observe the similarity between $S_n$ and the integrand in equation \eqref{equation:stat_specbreak} of the previous section: except for the difference regarding the calculation of pseudo-observations, the functionals only differ by the factor $\sqrt{\lfloor \zeta n \rfloor(n-\lfloor \zeta n \rfloor)}/n$ which assigns less weight to change point candidates close to the sample's boundaries. Define $Z_n:=\{1/n, \ldots, (n-1)/n\}.$ 
We consider three alternative test statistics which pick the most extreme realization within the set $Z_n$ of change point candidates: \index{functional! Kolmogorov-Smirnov} \index{functional! Kuiper} \index{functional! Cram\'{e}r-von Mises}
\begin{align}
\label{equation:ubCvM}
& T_n^{1}=\max_{\zeta \in Z_n} \int_{[0,1]^d} \Sb_n(\zeta,\mathbf{u})^2 d\widehat{C}_n(\mathbf{u}), \\  \label{equation:ubK}
& T_n^{2}=\max_{\zeta \in Z_n}  \left\{ \max_{\mathbf{u} \in \{\widehat{\mathbf{U}}_j\}_{j=1, \ldots,n}} \Sb_n(\zeta,\mathbf{u})- \min_{\mathbf{u} \in \{\widehat{\mathbf{U}}_j\}_{j=1, \ldots,n}} \Sb_n(\zeta,\mathbf{u}) \right\}, \\ \label{equation:ubKS}
& T_n^{3}=\max_{\zeta \in Z_n}  \left\{ \max_{\mathbf{u} \in \{\widehat{\mathbf{U}}_j\}_{j=1, \ldots,n}} | \Sb_n(\zeta,\mathbf{u}) | \right\},
\end{align}
which are the {\it maximally selected Cram\'{e}r-von Mises (CvM), Kuiper (K), and Kolmogorov-Smirnov (KS) statistic}\index{maximally selected statistic! Cram\'{e}r-von Mises}\index{maximally selected statistic! Kuiper}\index{maximally selected statistic! Kolmogorov-Smirnov}, respectively. We refer to \citet{horvathshao07} for an investigation of these statistics in a univariate context based on independent and identically distributed observations; $T_n^3$ is investigated in \citet{inoue01} for general multivariate distribution functions under strong mixing conditions as well as in \citet{remillard10} with an application to the copula of GARCH residuals.

\begin{theorem}
\label{theorem:kieferasymptotics}
Consider a sample $\mathbf{X}_{1},\ldots,\mathbf{X}_{n}$ of a strictly stationary process $(\mathbf{X}_j)_{j \in \mathbb{Z}}$ satisfying the strong mixing condition $\alpha_{\mathbf{X}}(r)=\mathcal{O}(r^{-4-d\{1+\eps\}})$ for some $0 < \eps \leq 1/4$.
Then under the null hypothesis, in $\left(\ell^\infty([0,1]^{d+1}),\|\cdot\|_\infty \right),$
\begin{align*}	
	\Sb_n(\zeta,\mathbf{u}) \overset{w.}{\longrightarrow} \Sb_C(\zeta, \mathbf{u})= \BB_C(\zeta,\mathbf{u})- \zeta \BB_C(1,\mathbf{u}),
\end{align*}
where $\BB_C(\zeta,\mathbf{u})$ denotes a (centered) $C$-Kiefer process,
with covariance structure
\begin{align*}
\Cov(\BB_C(\zeta_1, \mathbf{u}),\BB_C(\zeta_2, \mathbf{v})) = \min(\zeta_1,\zeta_2) \sum_{j \in \mathbb{Z}} \Cov\left(\mathbf{1}_{\{\mathbf{U}_0 \leq \mathbf{u}\}},\mathbf{1}_{\{\mathbf{U}_j \leq \mathbf{v} \}}\right)
\end{align*}
for all $\zeta_1,\zeta_2 \in [0,1]$ and $\mathbf{u},\mathbf{v} \in [0,1]^d.$ This in particular implies weak convergence of the test statistics $T_n^1,$ $T_n^2,$ and $T_n^3$ under $H_0$:
\begin{align*}
& T_n^1 \overset{w.}{\longrightarrow} \sup_{0 \leq \zeta \leq 1} \int_{[0,1]^d} \{\BB_C(\zeta,\mathbf{u})- \zeta \BB_C(1,\mathbf{u})\}^2 dC(\mathbf{u}), \\
& T_n^{2} \overset{w.}{\longrightarrow} \sup_{0 \leq \zeta \leq 1}  \left[ \sup_{\mathbf{u} \in [0,1]^d} \{\BB_C(\zeta,\mathbf{u})- \zeta \BB_C(1,\mathbf{u})\}- \inf_{\mathbf{u} \in [0,1]^d} \{\BB_C(\zeta,\mathbf{u})- \zeta \BB_C(1,\mathbf{u})\} \right], \\
& T_n^{3} \overset{w.}{\longrightarrow} \sup_{0 \leq \zeta \leq 1} \left[ \sup_{\mathbf{u} \in [0,1]^d} | \{\BB_C(\zeta,\mathbf{u})- \zeta \BB_C(1,\mathbf{u})\} | \right].
\end{align*}
\end{theorem}


%
	
Similar calculations as in \cite{horvathshao07} reveal that $T_n^i \rightarrow \infty$ for $i=1,2,3$ under $H_1$.
Hence, a test which rejects $H_0$ for unlikely large values of $T_n^i$ is consistent against general alternatives. Approximate critical values of the tests can be derived from the tapered block multiplier technique.

\begin{proposition}
 \label{proposition:estimationSunspec}
Consider a sample $\mathbf{X}_{1},\ldots,\mathbf{X}_{n}$ of a process $(\mathbf{X}_j)_{j \in \mathbb{Z}}$ which satisfies $X_{j,i} \sim F_i$ for all $j\in \mathbb{Z}$ and $i=1,\ldots, d.$ Further assume the process to fulfill the strong mixing assumptions of Theorem \ref{theorem:taperedmultiplier}. Let $\xi_{1,n}, \ldots, \xi_{n,n}$ denote a sample of a tapered block multiplier process $(\xi_{j,n})_{j \in \mathbb{Z}}$ satisfying \bf A1\it, \bf A2\it, \bf A3 \it with block length $l(n) \rightarrow \infty,$ where $l(n)=\mathcal{O}(n^{1/2-\eps})$ for $0< \eps <1/2$ and define
$
	\widehat{\Sb}_n(\zeta,\mathbf{u}) = \widehat{\BB}_n^M(\zeta,\mathbf{u}) - \zeta \widehat{\BB}_n^M(1,\mathbf{u}),
$
where
\begin{align*}
	\widehat{\BB}_n^M(\zeta,\mathbf{u}) = n^{-1/2} \sum_{j=1}^{\lfloor \zeta n\rfloor} \left( \frac{ \xi_{n,j}}{\bar \xi_{\lfloor \zeta n\rfloor}} - 1\right) \mathbf{1}_{\{ \widehat{\mathbf{U}}_j \le \mathbf{u} \}}.
\end{align*}
Then, under the null hypothesis, $\widehat{\Sb}_n \weakcondp{\xi} \Sb_C$, while under the alternative $\widehat{\Sb}_n=O_\PP(1)$.
\end{proposition}

An application of the continuous mapping theorem proves consistency of the tapered block multiplier-based tests. The p-values of the test statistics are estimated as shown in Equation (\ref{equation:estimate_p}).

For simplicity, the change point location is assessed under the assumption that there is at most one change point. In this case, the alternative hypothesis can as well be formulated:
\begin{align*}
& H_{1b}: \exists \lambda \in [0,1] \ \text{such that} \ \mathbf{U_j} \sim \begin{cases}   C_1 \ \text{for all} \ j=1, \ldots, \lfloor \lambda n \rfloor, \\ C_2 \ \text{for all} \ j=\lfloor \lambda n \rfloor +1, \ldots, n, \end{cases}
\end{align*}
where $C_1$ and $C_2$ are assumed to differ on a non-empty subset of $[0,1]^d.$ An estimator for the {\it location of the change point} $\widehat{\lambda}^i_n,$ $i=1, 2, 3,$ is obtained by replacing $\max$ functions by $\arg\max$ functions in Equations (\ref{equation:ubCvM}), (\ref{equation:ubK}), and (\ref{equation:ubKS}). For ease of exposition, the superindex $i$ is dropped in the following if no explicit reference to the functional is required. Given a (not necessarily correct) change-point estimator $\widehat{\lambda}_n,$ the empirical copula of $\mathbf{X}_1, \ldots, \mathbf{X}_{\lfloor \widehat{\lambda}_n n \rfloor}$ is an estimator of the unknown mixture distribution given by
\begin{align}
\label{equation:lambdahatmix1}
C_{\widehat{\lambda}_n,1}(\mathbf{u})=\mathbf{1}_{\left\{\widehat{\lambda}_n \leq \lambda \right\}}C_{1}(\mathbf{u})+\mathbf{1}_{\left\{\widehat{\lambda}_n > \lambda \right\}}
\widehat{\lambda}_n^{-1}
\left[\lambda C_{1}(\mathbf{u})+\left(\widehat{\lambda}_n-\lambda\right)C_{2}(\mathbf{u})\right]
\end{align}
\citep[for an analogous estimator related to general distribution functions, see][]{carlstein88}. The latter coincides with $C_1$ if and only if the change point is estimated correctly. On the other hand, the empirical copula of $\mathbf{X}_{\lfloor \widehat{\lambda}_n n \rfloor+1}, \ldots, \mathbf{X}_n$ is an estimator of the unknown mixture distribution given by
\begin{align}
\label{equation:lambdahatmix2}
C_{\widehat{\lambda}_n,2}(\mathbf{u})= & \mathbf{1}_{\left\{\widehat{\lambda}_n \leq \lambda \right\}}
(1-\widehat{\lambda}_n)^{-1} \left[(\lambda-\widehat{\lambda}_n) C_{1}(\mathbf{u})+\left(1-\lambda\right)
C_{2}(\mathbf{u})\right] \\ \nonumber
& +\mathbf{1}_{\left\{\widehat{\lambda}_n > \lambda \right\}} C_{2}(\mathbf{u}),
\end{align}
for all $\mathbf{u} \in [0,1]^d.$ The latter coincides with $C_2$ if and only if the change point is estimated correctly. Consistency of $\widehat{\lambda}_n$ follows from consistency of the empirical copula and the fact that the difference of the two mixture distributions given in Equations (\ref{equation:lambdahatmix1}) and (\ref{equation:lambdahatmix2}) is maximal in the case $\widehat{\lambda}_n=\lambda.$
\citet{bai97} iteratively applies the setting considered above to test for multiple breaks (one at a time), indicating a direction of future research to estimate locations of multiple change points in the dependence structure.

\bigskip

\bf Finite sample properties. \rm Size and power of the tests for a constant copula are shown in Tables \ref{tab:powerubn400} and \ref{tab:powerubn800} for $n=400$ and $n=800$, respectively. The observations are either serially independent or from a strictly stationary AR($1$) process where the univariate innovations are either linked by a Clayton or a Gumbel--Hougaard copula. We consider the alternative hypothesis $H_{1b}$ of at most one unspecified change point. If present, then this change point is located after observation $\lfloor \lambda n \rfloor=n/2$ and only affects the parameter within the investigated Clayton or Gumbel--Hougaard families: the
copula $C_1$ is parameterized such that Kendall's $\tau_1=0.2$, the copula $C_2$ such that Kendall's $\tau_2 \in \{0.2, 0.6, 0.9\}$. We consider $S=1,000$ ($n=400$) or $S=500$  ($n=800$) tapered block multiplier simulations based on normal multiplier random variables with block length $l_M(400)=5$, $l_M(800)=6$ and  kernel function $\kappa_2$.

\begin{table}[htbp!]
\caption{{\small{\bf Size and power of tests for a constant copula with unspecified change point candidate}. Results are based on $1,000$ Monte Carlo replications, $n=400,$ $S=1,000,$ kernel function $\kappa_2,$ and $\alpha=5\%$; additionally, the estimated change point location $\widehat{\lambda}_n,$ $\hat{\sigma}(\widehat{\lambda}_n),$ and $MSE(\widehat{\lambda}_n)\times 10^2$ are reported.}}
\label{tab:powerubn400}
\begin{center}
\small
\begin{tabular}{lllrrr|rr|rr|rr} \hline
 & & & & & & & & & \\
& & & \multicolumn{3}{l|}{\bf size/power \rm} & \multicolumn{2}{l|}{\bf $\widehat{\lambda}_n$ \rm} & \multicolumn{2}{l|}{\bf $\hat{\sigma}\left(\widehat{\lambda}_n\right)$ \rm}
& \multicolumn{2}{l}{\bf $MSE\left(\widehat{\lambda}_n\right)$ \rm} \\
& & & & & & & & &  \vspace{-10pt} \\
 & & \multicolumn{1}{r}{$\tau_2$} & $0.2$ & $0.6$ & $0.9$ &  $0.6$ & $0.9$ & $0.6$ & $0.9$ & $0.6$ & $0.9$ \\
\hline \hline
& & & & & & & & & & &  \\
\multicolumn{5}{l}{\bf i.i.d. setting \rm} & & & & & & & \\ 
Clayton & $l=1$ & $CvM$ & 0.061 & 0.406 & 0.873 & 0.511 & 0.506 & 0.093 & 0.061 & 0.871 & 0.372 \\
& & $K$ &  0.040 & 0.495 & 0.991 & 0.496 & 0.496 & 0.074 & 0.048 & 0.556 & 0.228 \\
& & $KS$ & 0.062 & 0.409 & 0.905 & 0.507 & 0.503 & 0.079 & 0.056 & 0.627 & 0.319 \\
& $l=5$ & $CvM$ & 0.035 & 0.342 & 0.847 & 0.519 & 0.509 & 0.084 & 0.058 & 0.735 & 0.349 \\
& & $K$ & 0.031 & 0.375 & 0.983 & 0.495 & 0.496 & 0.070 & 0.049 & 0.489 & 0.246 \\
& & $KS$ & 0.036 & 0.337 & 0.881 & 0.506 & 0.504 & 0.080 & 0.055 & 0.645 & 0.306 \\
Gumbel & $l=1$ & $CvM$ & 0.047 & 0.456 & 0.903 & 0.509 & 0.507 & 0.083 & 0.055 & 0.694 & 0.304 \\
& & $K$ & 0.056 & 0.474 & 0.992 & 0.494 & 0.495 & 0.074 & 0.050 & 0.548 & 0.251 \\
& & $KS$ & 0.052 & 0.437 & 0.910 & 0.502 & 0.504 & 0.080 & 0.054 & 0.644 & 0.291 \\
& $l=5$ & $CvM$ & 0.043 & 0.387 & 0.880 & 0.506 & 0.507 & 0.086 & 0.055 & 0.739 & 0.310 \\
& & $K$ & 0.026 & 0.343 & 0.974 & 0.487 & 0.496 & 0.072 & 0.046 & 0.518 & 0.217 \\
& & $KS$ & 0.041 & 0.349 & 0.884 & 0.497 & 0.503 & 0.079 & 0.056 & 0.642 & 0.312 \\
& & & & & & & & & & & \\

\multicolumn{5}{l}{\bf AR($1$) setting with $\beta=0.25$ \rm} & & & & & & & \\ 
Clayton & $l=1$ & $CvM$ & 0.187 & 0.487 & 0.846 & 0.519 & 0.510 & 0.116 & 0.075 & 1.390 & 0.566  \\
& & $K$ & 0.139 & 0.562 & 0.994 & 0.490 & 0.494 & 0.084 & 0.052 & 0.707 & 0.276 \\
& & $KS$ & 0.177 & 0.499 & 0.913 & 0.504 & 0.506 & 0.098 & 0.071 & 0.969 & 0.503  \\
& $l=5$ & $CvM$ & 0.046 & 0.243 & 0.697 & 0.516 & 0.512 & 0.102 & 0.070 & 1.073 & 0.497  \\
& & $K$ & 0.039 & 0.289 & 0.954 & 0.496 & 0.493 & 0.073 & 0.056 & 0.535 & 0.315  \\
& & $KS$ & 0.040 & 0.253 & 0.771 & 0.506 & 0.506 & 0.095 & 0.065 & 0.901 & 0.421  \\
Gumbel & $l=1$ & $CvM$ & 0.185 & 0.536 & 0.878 & 0.515 & 0.513 & 0.119 & 0.079 & 1.443 & 0.649 \\
& & $K$ & 0.123 & 0.576 & 0.997 & 0.487 & 0.493 & 0.086 & 0.051 & 0.746 & 0.272 \\
& & $KS$ & 0.167 & 0.541 & 0.913 & 0.504 & 0.509 & 0.103 & 0.075 & 1.078 & 0.573 \\
& $l=5$ & $CvM$ & 0.042 & 0.295 & 0.706 & 0.519 & 0.506 & 0.096 & 0.074 & 0.949 & 0.547 \\
& & $K$ & 0.040 & 0.294 & 0.939 & 0.495 & 0.495 & 0.084 & 0.053 & 0.716 & 0.282 \\
& & $KS$ & 0.050 & 0.287 & 0.745 & 0.509 & 0.504 & 0.089 & 0.067 & 0.793 & 0.457 \\
& & & & & & & & & & & \\
\hline
\end{tabular}
\normalsize
\end{center}
\end{table}

\begin{table}[htbp!]
\caption{{\small{\bf Size and power of tests for a constant copula with unspecified change point candidate}. Results are based on $1,000$ Monte Carlo replications, $n=800,$ $S=500,$ kernel function $\kappa_2,$ and $\alpha=5\%$; additionally, the estimated change point location $\widehat{\lambda}_n,$ $\hat{\sigma}(\widehat{\lambda}_n),$ and $MSE(\widehat{\lambda}_n)\times 10^2$ are reported.}}
\label{tab:powerubn800}
\begin{center}
\small
\begin{tabular}{lllrrr|rr|rr|rr} \hline
 & & & & & & & & & \\
& & & \multicolumn{3}{l|}{\bf size/power \rm} & \multicolumn{2}{l|}{\bf $\widehat{\lambda}_n$ \rm} & \multicolumn{2}{l|}{\bf $\hat{\sigma}\left(\widehat{\lambda}_n\right)$ \rm}
& \multicolumn{2}{l}{\bf $MSE\left(\widehat{\lambda}_n\right)$ \rm} \\
& & & & & & & & &  \vspace{-10pt} \\
 & & \multicolumn{1}{r}{$\tau_2$} & $0.2$ & $0.6$ & $0.9$ &  $0.6$ & $0.9$ & $0.6$ & $0.9$ & $0.6$ & $0.9$ \\
\hline \hline
& & & & & & & & & & &  \\
\multicolumn{5}{l}{\bf i.i.d. setting \rm} & & & & & & & \\ 
Clayton & $l=1$ & $CvM$ & 0.033 & 0.695 & 0.999 & 0.508 & 0.506 & 0.079 & 0.041 & 0.623 & 0.173 \\
& & $K$ &  0.054 & 0.901 & 1.000 & 0.494 & 0.497 & 0.059 & 0.029 & 0.350 & 0.087 \\
& & $KS$ & 0.046 & 0.734 & 0.999 & 0.505 & 0.503 & 0.068 & 0.038 & 0.463 & 0.142 \\
& $l=6$ & $CvM$ & 0.043 & 0.681 & 0.997 & 0.510 & 0.505 & 0.076 & 0.042 & 0.585 & 0.177 \\
& & $K$ & 0.037 & 0.838 & 1.000 & 0.494 & 0.497 & 0.057 & 0.031 & 0.335 & 0.096 \\
& & $KS$ & 0.036 & 0.699 & 0.998 & 0.507 & 0.503 & 0.062 & 0.039 & 0.384 & 0.149 \\
Gumbel & $l=1$ & $CvM$ & 0.049 & 0.510 & 0.953 & 0.510 & 0.506 & 0.090 & 0.049 & 0.823 & 0.243 \\
& & $K$ & 0.037 & 0.688 & 1.000 & 0.490 & 0.496 & 0.061 & 0.032 & 0.372 & 0.106 \\
& & $KS$ & 0.051 & 0.509 & 0.966 & 0.502 & 0.503 & 0.081 & 0.044 & 0.672 & 0.199 \\
& $l=6$ & $CvM$ & 0.055 & 0.721 & 0.998 & 0.505 & 0.504 & 0.069 & 0.037 & 0.477 & 0.138 \\
& & $K$ & 0.034 & 0.775 & 1.000 & 0.498 & 0.496 & 0.059 & 0.029 & 0.355 & 0.086 \\
& & $KS$ & 0.045 & 0.686 & 0.998 & 0.504 & 0.500 & 0.065 & 0.039 & 0.428 & 0.156 \\
& & & & & & & & & & & \\

\multicolumn{5}{l}{\bf AR($1$) setting with $\beta=0.25$ \rm} & & & & & & &  \\ 
Clayton & $l=1$ & $CvM$ & 0.184 & 0.712 & 0.995 & 0.519 & 0.507 & 0.097 & 0.056 & 0.973 & 0.321  \\
& & $K$ & 0.122 & 0.914 & 1.000 & 0.495 & 0.497 & 0.063 & 0.034 & 0.411 & 0.122 \\
& & $KS$ & 0.169 & 0.756 & 0.999 & 0.512 & 0.506 & 0.084 & 0.050 & 0.714 & 0.248  \\
& $l=6$ & $CvM$ & 0.065 & 0.488 & 0.939 & 0.506 & 0.508 & 0.089 & 0.054 & 0.803 & 0.293 \\
& & $K$ & 0.057 & 0.724 & 1.000 & 0.493 & 0.496 & 0.060 & 0.033 & 0.359 & 0.113 \\
& & $KS$ & 0.062 & 0.521 & 0.965 & 0.500 & 0.503 & 0.072 & 0.045 & 0.527 & 0.206 \\
Gumbel & $l=1$ & $CvM$ & 0.207 & 0.760 & 0.992 & 0.507 & 0.509 & 0.091 & 0.052 & 0.841 & 0.273 \\
& & $K$ & 0.141 & 0.879 & 1.000 & 0.489 & 0.497 & 0.070 & 0.036 & 0.496 & 0.135 \\
& & $KS$ & 0.182 & 0.780 & 0.998 & 0.500 & 0.505 & 0.089 & 0.047 & 0.805 & 0.218 \\
& $l=6$ & $CvM$ & 0.052 & 0.533 & 0.959 & 0.514 & 0.508 & 0.083 & 0.052 & 0.707 & 0.273 \\
& & $K$ & 0.046 & 0.670 & 1.000 & 0.495 & 0.496 & 0.065 & 0.032 & 0.426 & 0.106 \\
& & $KS$ & 0.053 & 0.520 & 0.974 & 0.508 & 0.505 & 0.076 & 0.048 & 0.584 & 0.232 \\
& & & & & & & & & & & \\
\hline
\end{tabular}
\normalsize
\vspace{-12pt}
\end{center}
\end{table}

In the case of i.i.d. observations, we observe that the tapered block multiplier works similarly well as the standard multiplier (i.e., $l_M=1$): the asymptotic size of the test, chosen to be $5\%,$ is well approximated and its power increases in the difference $\tau_2-\tau_1.$ The estimated location of the change point, $\widehat{\lambda}_n,$ is close to its theoretical value. Moreover, its standard deviation $\hat{\sigma}(\widehat{\lambda}_n)$ as well as its mean squared error MSE($\widehat{\lambda}_n$) are decreasing in the difference $\tau_2-\tau_1.$ In the case of serially dependent observations sampled from AR($1$) processes with $\beta=0.25,$ we find that the observed size of the test strongly deviates from its nominal size (chosen to be $5\%$) if serial dependence
is neglected and the block length $l_M=1$ is used: its estimates are reaching up to $18.7\%.$ The test based on the tapered block multiplier with block length $l_M(400)=5$ yields rejection quotas which approximate the asymptotic size well in all settings considered. These results are strengthened in Table \ref{tab:powerubn800} which shows results of MC simulations for sample size $n=800$ and block length $l_M=6.$
The power improves considerably with the increased amount of observations and the change point location is well captured. Standard deviation and mean squared error of the estimated location of the change point, $\widehat{\lambda}_n,$ decrease in the difference $\tau_2-\tau_1.$

Comparing the tests based on statistics $T_n^1,$ $T_n^2,$ and $T_n^3,$ we find that the test based on the Kuiper-type statistic performs best. The results indicate that the nominal size is well approximated in finite samples and that the test is most powerful in many settings. Likewise, with regard to the estimated location of the change point, the Kuiper-type statistic performs best in mean and in mean squared error.

The introduced tests for a constant copula offer some connecting factors for further research. For instance, \citet{inoue01} investigates nonparametric change point tests for the joint distribution of strongly mixing random vectors and finds that the observed size of the test heavily depends on the choice of the block length $l$ in the resampling procedure. For different types of serially dependent observations, e.g., AR($1$) processes with higher coefficient for the lagged variable or GARCH($1,1$) processes, it is of interest to investigate the optimal choice of the block length for the tapered block multiplier-based test with unspecified change point candidate. Moreover, test statistics based on different functionals offer potential for improvements. For instance, the Cram\'{e}r-von Mises functional introduced by \citet{remillardscaillet09} led to strong results in the case of a specified change point candidate. Though challenging from a computational point of view, an application of this functional to the case of unspecified change point candidate(s) is of interest as the functional yields very powerful tests.

\section{Conclusion}\label{section:conclusion}

Consistent tests for constancy of the copula with specified or unspecified change point candidate are introduced. We observe a trade-off in assumptions required for the testing: if a change point candidate is specified, then the test is consistent whether or not there is a simultaneous change point in marginal distribution function(s). If change point candidate(s) are unspecified, then the assumption of strictly stationary marginal distribution functions is required and allows to drop continuity assumptions on the partial derivatives of the underlying copula(s). Tests are shown to behave well in size and power when applied to various types of dependent observations. P-Values of the tests are estimated using a tapered block multiplier technique which is based on serially dependent multiplier random variables; the latter is shown to perform better than the block bootstrap in mean and mean squared error when estimating the asymptotic covariance structure of the empirical copula process in various settings.

\textbf{Acknowledgements.} The authors would like to thank two unknown referees and the Associate editor for their constructive comments on an earlier version of this manuscript. We are indebted to Friedrich Schmid and Ivan Kojadinovic for helpful comments and discussions. Moreover, we are grateful to the conference participants at the ``German Open Conference on Probability and Statistics 2010'' (Leipzig, Germany) and at the ``7th Conference on Multivariate Distributions with Applications 2010'' (Maresias, Brazil) for constructive comments on an earlier version of this paper. We would like to thank the Regional Computing Center at the University of Cologne for providing the computational resources required. Martin Ruppert gratefully acknowledges financial support by the German Research Foundation (DFG). Axel B\"ucher is thankful for financial support through the collaborative research center ``Statistical modeling of nonlinear dynamic processes'' (SFB 823) of the German Research Foundation (DFG).


\appendix
\section{Proofs}
\label{appendix:proofTheorems}

\begin{proof}[{\it Proof of Theorem \ref{theorem:empcopprocess}}]
Since $\widehat{F}_i(X_{j,i})\le u_i$ if and only if $\widehat{G}_i(U_{j,i})\le u_i$, where $\widehat{G}_i$ denotes the empirical distribution function of $U_{1,i},\dots,U_{n,i}$, we can assume without loss of generality, that $\mathbf{X}_j=\mathbf{U}_j\sim C$.
Weak convergence of $\BB_{C,n}$ under the strong mixing condition $\alpha_\mathbf{X}(r)=\mathcal{O}(r^{-a})$ for some $a>1$  is established in \citet{rio00}.
Thus, Condition 2.1 in \citet{buecvolg2011} is satisfied and an application of the functional delta method and of Theorem 2.4 in \cite{buecvolg2011} yields
\[
	\sqrt{n}\{C_n(C_{n,1}^{-1}, \dots, C_{n,d}^{-1}) - C\} \weak \GG_C,
\]
where $C_{n,i}$ denotes the $i$-th marginal of $C_n$, $i=1,\dots,d$.
The assertion follows from $\|\widehat{\GG}_{C,n}-\sqrt{n}\{C_n(C_{n,1}^{-1}, \dots, C_{n,d}^{-1}) - C\}\|_\infty = O(n^{-1/2})$.
\end{proof}

\begin{proof}[{\it Proof of Theorem \ref{theorem:blockbootstrap}}]
Without loss of generality we may assume $\mathbf{X}_j=\mathbf{U}_j$.
Let $C_{n,b}$ denote the empirical distribution function of the sample
\begin{align*}
\mathbf{U}_{H_1+1},\dots,\mathbf{U}_{H_1+l_B},\mathbf{U}_{H_2+1},\dots,\mathbf{U}_{H_2+l_B}, \dots \dots, \mathbf{U}_{H_k+1}, \dots, \mathbf{U}_{H_k+l_B}.
\end{align*}
Again, $\|\widehat{\GG}_{C,n}^B-\sqrt{n}\{C_{n,b}(C_{n,b,1}^{-1}, \dots, C_{n,b,d}^{-1}) - C_n\}\|_\infty=O(n^{-1/2})$, and the
result follows from Corollary 2.11 in \cite{buecvolg2011}.
\end{proof}

\begin{proof}[{\it Proof of Theorem \ref{theorem:taperedmultiplier}}]
Without loss of generality we may assume $\mathbf{X}_j=\mathbf{U}_j$.
The process
\[
	\BB_{C,n}^M(\mathbf{u})=\sqrt{n} \left\{ \frac{1}{n} \sum_{j=1}^{n} \frac{\xi_{n,j}}{\bar{\xi}_n} \mathbf{1}_{\{\mathbf{U}_{j} \leq \mathbf{u}\}}-C_{n}(\mathbf{u})\right\}
		 = \frac{1}{\sqrt{n}} \sum_{j=1}^{n} \left( \frac{\xi_{n,j}}{\bar{\xi}_n} -1 \right) \left\{\mathbf{1}_{\{\mathbf{U}_{j} \leq \mathbf{u}\}} - C(\mathbf{u}) \right\}
\]
is defined on a product space $(\Omega_1\times\Omega_2,\Fc_1\otimes\Fc_2,\PP_1\otimes\PP_2)$, where $\mathbf{X}_j$ is defined on $(\Omega_1,\Fc_1,\PP_1)$ and $\xi_{n,j}$ on $(\Omega_2,\Fc_2,\PP_2)$ for all $j=1,\dots,n$.
It follows from the results in \citet[Section 3.3]{buehlmann93} and \cite[Section 1.5]{vaartwellner96} that there exists a set $A_1\in\Fc_1$ with $\PP_1(A_1)=1$ such that
\[
	\BB_{C,n}^M(\omega_1,\cdot) \weak \BB_C \text{ in } \ell^\infty([0,1]^d) \quad\text{for all } \omega_1\in A_1,
\]
where $\weak$ denotes weak convergence with respect to the multipliers.
By Theorem 1.5.7 and its addendum in \cite{vaartwellner96}, $\BB_{C,n}^M(\omega_1,\cdot)$ is asymptotically uniformly equicontinuous for each such $\omega_1$. As $\sup_{u_i\in[0,1]}|\widehat{F}_{i}^-(u_i) - u_i| =o(1)$, $\PP_1$-almost surely, we obtain
\[
	\tilde\BB_{C,n}^M(\mathbf{u}) = \BB_{C,n}^M(\widehat F_1^{-1}(u_1), \dots, \widehat F_d^{-1}(u_d)) \weak \BB_{C} \quad (\PP_1\text{-a.s.}).
\]
Note that
\[
	\|\tilde\BB_{C,n}^M - \widehat{\BB}_{C,n}^M\|_\infty =O( n^{-1/2}\max_{j=1}^n |\xi_{n,j}/\bar \xi_n|) = o_{\PP_2}(1),
\]
where the last estimation follows from $\PP_2(|\xi_{n,1}|>x)=o(x^{-2})$ and $\PP_2(\max_{j=1}^n |\xi_{n,j}| >\eps \sqrt{n} \bar \xi_n) \le \PP_2(\bar \xi _n \le 1/2) +
n\PP_2(|\xi_{n,1}|>\eps \sqrt{n}/2).$ Thus $\widehat{\BB}_{C,n}^M  \weak \BB_{C}$, $\PP_1$-almost surely,
and since $\widehat{\BB}_{C,n}^M$ is a ball-measurable element of the space of cadlag functions $\left(D([0,1]^{d}),\|\cdot\|_\infty \right)\subset(\ell^\infty([0,1]^d),\|\cdot\|_\infty)$, a detailed look at Theorem 1.7.2 and the proof of Theorem 2.9.6 in \cite{vaartwellner96} allows to translate this result into the Hoffmann-J\o rgensen weak convergence, i.e.,
\[
	\widehat{\BB}_{C,n}^M \weakcondp{\xi} \BB_C \quad\text{in }(\ell^\infty([0,1]^d),\|\cdot\|_\infty).
\]
Finally, set
\[
	\tilde{\GG}_{C,n}^M (\mathbf{u}) = \widehat{\BB}_{C,n}^M(\mathbf{u})  - \sum_{i=1}^d D_iC(\mathbf{u}) \widehat{\BB}_{C,n}^M(\mathbf{u}^{(i)}).
\]
Proceeding as in the proof of Proposition 3.2 in \cite{segers11} we obtain $\| \tilde{\GG}_{C,n}^M - \widehat{\GG}_{C,n}^M \|_\infty = o_P(1)$, which yields the assertion of the Theorem after an application of Lemma B.1 in \cite{buedetvol2011}.
\end{proof}


\begin{proof}[{\it Proof of Theorem \ref{theorem:asymptoticsS}}]
Without loss of generality we may assume $\mathbf{X}_j=\mathbf{U}_j$. We begin by proving the joint weak convergence
\begin{align} \label{eq:jointconv}
	(\BB_{C,n,1},\BB_{C,n,2}) \weak (\BB_{C,1}, \BB_{C,2})
\end{align}
in $\ell^\infty([0,1]^d) \times \ell^\infty([0,1]^d)$, where
\begin{align*}
	\BB_{C,n,1}(\mathbf{u}) &= \frac{1}{\sqrt{\lfloor \lambda n \rfloor}} \sum_{j=1}^{\lfloor \lambda n \rfloor} \left\{ \mathbf{1}_{\{\mathbf{U}_{j} \leq \mathbf{u}\}}-C(\mathbf{u}) \right\} \\ 
	\BB_{C,n,2}(\mathbf{u}) &= \frac{1}{\sqrt{n - \lfloor \lambda n \rfloor}} \sum_{j=\lfloor \lambda n \rfloor + 1}^n \left\{ \mathbf{1}_{\{\mathbf{U}_{j} \leq \mathbf{u}\}}-C(\mathbf{u}) \right\}. 
\end{align*}
Asymptotic tightness of $(\BB_{C,n,1},\BB_{C,n,2})$ follows from Lemma 1.4.3 in \cite{vaartwellner96} and the fact that
\[
	\BB_{C,n,2} = \sqrt{\frac{1}{1-s}} \BB_{C,n} - \sqrt{\frac{s}{1-s}} \BB_{C,n,1} + o_\PP(1),
\]
where both processes on the right-hand side are asymptotically tight by Theorem \ref{theorem:empcopprocess}.
It therefore suffices to show that
\begin{multline*}
	\big(\BB_{C,n,1}(\mathbf{u_1}), \dots, \BB_{C,n,1}(\mathbf{u_k}), \BB_{C,n,2}(\mathbf{v_1}), \dots, \BB_{C,n,2}(\mathbf{v_l})\big) \\
		\weak \big(\BB_{C,1}(\mathbf{u_1}), \dots, \BB_{C,1}(\mathbf{u_k}), \BB_{C,2}(\mathbf{v_1}), \dots, \BB_{C,2}(\mathbf{v_l})\big)
\end{multline*}
in $\R^{k+l}$, see Problem 1.5.3 in \cite{vaartwellner96}. For the sake of a clear exposition we only consider the case $k=l=1$, the general case follows along similar lines. We have to prove that, for each $c=(c_1,c_2)\in \R^2$,
\[
	c_1\BB_{C,n,1}(\mathbf{u}) + c_2 \BB_{C,n,2}(\mathbf{v}) \weak c_1\BB_{C,1}(\mathbf{u}) + c_2 \BB_{C,2}(\mathbf{v}).
\]
Setting $Z_j(\mathbf{u})= \mathbf{1}_{\{\mathbf{U}_{j} \leq \mathbf{u}\}}-C(\mathbf{u})$, the left-hand side of the previous expression can be written as $\sum_{j=1}^n \zeta_{n,j}$, where
\[
	\zeta_{n,j} =\frac{c_1}{\sqrt{\lfloor \lambda n \rfloor}} Z_j(\mathbf{u})  \mathbf{1}_{\{j \le \lfloor \lambda n\rfloor\} } + \frac{c_2}{\sqrt{n-\lfloor \lambda n \rfloor}} Z_j(\mathbf{v}) \mathbf{1}_{\{j > \lfloor \lambda n\rfloor \}}.
\]
The asserted weak convergence follows from Theorem 2.1 in \cite{peligrad1996}, if we prove the corresponding conditions of that Theorem. We have
\begin{multline*}
	\EE\left[ \big(\sum_{j=1}^n \zeta_{n,j} \big)^2\right]
	= \frac{c_1^2}{\lfloor \lambda n \rfloor} \sum_{i,j=1}^{ \lfloor \lambda n \rfloor}  \EE[Z_i(\mathbf{u})Z_j(\mathbf{u})] +  \frac{c_2^2}{n-\lfloor \lambda n \rfloor} \sum_{i,j=\lfloor \lambda n \rfloor+1}^{n} \EE[Z_i(\mathbf{v})Z_j(\mathbf{v})]  \\
		+  \frac{2c_1c_2}{\sqrt{\lfloor \lambda n \rfloor(n-\lfloor \lambda n \rfloor)}} \sum_{i=1}^{\lfloor \lambda n \rfloor} \sum_{j=\lfloor \lambda n \rfloor+1}^{n} \EE[Z_i(\mathbf{u}) Z_j(\mathbf{v})].
\end{multline*}
The first two sums converge to $c_1^2\gamma(\mathbf{u},\mathbf{u})$ and $c_2^2 \gamma(\mathbf{v},\mathbf{v})$, respectively. Suppose that $\lfloor \lambda n \rfloor>n-\lfloor \lambda n \rfloor$, the opposite case is treated analogously. Then,
a tedious calculation shows that the third sum in the previous identity equals
\begin{align*}
	 \frac{2c_1c_2}{\sqrt{\lfloor \lambda n \rfloor ( n-\lfloor \lambda n \rfloor) }}
	\left\{ \sum_{k=1}^{n-\lfloor \lambda n \rfloor} k r(k) + \sum_{k=n-\lfloor \lambda n \rfloor+1}^{\lfloor \lambda n \rfloor-1} (n-\lfloor \lambda n \rfloor) r(k) + \sum_{k=\lfloor \lambda n \rfloor}^{n-1}(n-k)r(k) \right\},
\end{align*}
where $r(|i-j|)=\EE[Z_i(\mathbf{u})Z_j(\mathbf{v})]$.
The first sum in the curly bracket multiplied with $n^{-1}$ converges to $0$ by dominated convergence. If we multiply the other two sums with $n^{-1}$, then we obtain tails of absolutely converging series, which also converge to $0$. To conclude,
\[
	\EE\left[(\sum_{j=1}^n \zeta_{n,j})^2\right] \rightarrow  c^T \left(\begin{array}{c c} \gamma(\mathbf{u},\mathbf{u}) & 0  \\ 0  & \gamma(\mathbf{v},\mathbf{v})  \end{array} \right) c =\Var(c_1\BB_{C,1}(\mathbf{u}) + c_2 \BB_{C}(\mathbf{v})).
\]
Since $|\zeta_{n,j} | =O(n^{-1/2})$ uniformly in $j$, some easy calculations show that condition (2.1) and (2.2) in \cite{peligrad1996} are satisfied. Hence, by Theorem 2.1 in that reference, the asserted weak convergence follows.

Define
\begin{align*}
	\widehat{\GG}_{C,n,1}(\mathbf{u}) & = \frac{1}{\sqrt{\lfloor \lambda n \rfloor}} \sum_{j=1}^{\lfloor \lambda n \rfloor} \left\{ \mathbf{1}_{\{\widehat{\mathbf{U}}_{j} \leq \mathbf{u}\}}-C(\mathbf{u}) \right\} \\
	\widehat{\GG}_{C,n,2}(\mathbf{u}) & = \frac{1}{\sqrt{n - \lfloor \lambda n \rfloor}} \sum_{j=\lfloor \lambda n \rfloor + 1}^n \left\{ \mathbf{1}_{\{\widehat{\mathbf{V}}_{j} \leq \mathbf{u}\}}-C(\mathbf{u}) \right\}.
\end{align*}
An application of the functional delta method to $\BB_{C,n,p}$ with the mapping $\Phi$ from Theorem 2.4 in \cite{buecvolg2011} (and the usual estimation of the remainder term descending from the fact that $F_{nj}(X_{ij}) \le u_j  $ and $X_{ij} \le F_{nj}^{-1}(u_j)$ are not equivalent) yields
\[
	(\widehat{\GG}_{C,n,1},\widehat{\GG}_{C,n,2}) \weak (\GG_{C_1},\GG_{C_2}) \quad \text{in } \ell^\infty([0,1]^d)\times \ell^\infty([0,1]^d).
\]
Hence, Theorem \ref{theorem:asymptoticsS} follows from the continuous mapping theorem observing that, under $H_0$,
\begin{multline*}
	\sqrt{\frac{\lfloor \lambda n \rfloor (n-\lfloor \lambda n \rfloor )}{n}}
	\left\{\widehat{C}_{1,\dots, \lfloor \lambda n \rfloor }-\widehat{C}_{n-\lfloor \lambda n \rfloor, \dots, n} \right\}
		= \sqrt{\frac{n-\lfloor \lambda n \rfloor }{n}}\widehat{\GG}_{C,n,1}-\sqrt{\frac{\lfloor \lambda n \rfloor }{n}}\widehat{\GG}_{C,n,2}.
\end{multline*}
\end{proof}


\begin{proof}[{\it Proof of Proposition \ref{prop:estimationS}}]
Without loss of generality we may assume $\mathbf{X}_j=\mathbf{U}_j$. We consider the processes
\begin{align*}
	\BB_{C,n,1}^M(\mathbf{u}) &= \frac{1}{\sqrt{\lfloor \lambda n \rfloor}} \sum_{j=1}^{\lfloor \lambda n \rfloor} \left(\frac{\xi_{n,j}}{\bar \xi_n^{(1)} } - 1\right) \mathbf{1}_{\{\mathbf{U}_{j} \leq \mathbf{u}\}}   \\
	&= \frac{1}{\sqrt{\lfloor \lambda n \rfloor}} \sum_{j=1}^{\lfloor \lambda n \rfloor} \left(\frac{\xi_{n,j}}{\bar \xi_n^{(1)} } - 1\right) \left\{ \mathbf{1}_{\{\mathbf{U}_{j} \leq \mathbf{u}\}}  - C(\mathbf{u}) \right\} \\
	\BB_{C,n,2}^M(\mathbf{u}) &= \frac{1}{\sqrt{n - \lfloor \lambda n \rfloor}} \sum_{j=\lfloor \lambda n \rfloor+1}^n  \left(\frac{\xi_{n,j}}{\bar \xi_n^{(2)} } - 1\right)  \mathbf{1}_{\{\mathbf{U}_{j} \leq \mathbf{u}\}} \\
	&= \frac{1}{\sqrt{n - \lfloor \lambda n \rfloor}} \sum_{j=\lfloor \lambda n \rfloor+1}^n  \left(\frac{\xi_{n,j}}{\bar \xi_n^{(2)} } - 1\right) \left\{ \mathbf{1}_{\{\mathbf{U}_{j} \leq \mathbf{u}\}}  - C(\mathbf{u}) \right\}
\end{align*}
which are defined on the product space $(\Omega_1\times\Omega_2,\Fc_1\otimes\Fc_2,\PP_1\otimes\PP_2)$, where $\mathbf{X}_j$ is defined on $(\Omega_1,\Fc_1,\PP_1)$ and $\xi_{n,j}$ on $(\Omega_2,\Fc_2,\PP_2)$ for all $j=1,\dots,n$. For the subsequent investigations, we may replace the arithmetic means $\bar \xi_n^{(p)}$ by $1=\EE \xi_{n,j}$. As in the proof of Theorem \ref{theorem:taperedmultiplier}, the results in
\citet[Section 3.3]{buehlmann93} guarantee the existence of a set $A_1\in\Fc_1$ with $\PP_1(A_1)=1$ such that both
\[
	\BB_{C,n,1}^M(\omega_1,\cdot) \weak \BB_{C,1}, \qquad \BB_{C,n,2}^M(\omega_1,\cdot) \weak \BB_{C,2} \quad\text{for all } \omega_1\in A_1.
\]
(use an index shift for $\BB_{C,n,2}$), where $\weak$ denotes weak convergence with respect to the multipliers. In the following we show joint weak convergence with respects to the multipliers, $\PP_1$-almost surely. Joint asymptotic tightness follows from asymptotic tightness of the components. As in the proof of Theorem \ref{theorem:asymptoticsS} it remains to show
\begin{multline*}
	\big(\BB_{C,n,1}^M(\mathbf{u_1}), \dots, \BB_{C,n,1}^M(\mathbf{u_k}), \BB_{C,n,2}^M(\mathbf{v_1}), \dots, \BB_{C,n,2}^M(\mathbf{v_l})\big)(\omega_1, \cdot) \\
		\weak \big(\BB_{C,1}(\mathbf{u_1}), \dots, \BB_{C,1}(\mathbf{u_k}), \BB_{C,2}(\mathbf{v_1}), \dots, \BB_{C,2}(\mathbf{v_l})\big)
\end{multline*}
in $\R^{k+l}$, which can be done by Theorem 2.1 in \cite{peligrad1996}. Again, for the sake of a clear exposition, we only consider the case $k=l=1$. We have to prove that, for each $c=(c_1,c_2)\in \R^2$,
\[
	c_1\BB_{C,n,1}^M(\mathbf{u}) + c_2 \BB_{C,n,2}^M(\mathbf{v}) \weak c_1\BB_{C,1}(\mathbf{u}) + c_2 \BB_{C,2}(\mathbf{v}) \quad (\PP_1\text{-a.s.}).
\]
The left-hand side of this display can be written as $\sum_{j=1}^n \zeta_{n,j}$, where
\[
	\zeta_{n,j}
		= \frac{c_1}{\sqrt{\lfloor \lambda n \rfloor}} \mathbf{1}_{\{ \mathbf{U}_j\le \mathbf{u} \}}  \mathbf{1}_{\{j \le \lfloor \lambda n\rfloor\} }
		+ \frac{c_2}{\sqrt{n-\lfloor \lambda n \rfloor}} \mathbf{1}_{\{ \mathbf{U}_j\le \mathbf{v} \}} \mathbf{1}_{\{j > \lfloor \lambda n\rfloor \}}.
\]
An easy calculation shows that, uniformly in $\mathbf{u},\mathbf{v}$ and $\omega_1$,
\begin{align*} 	
	\left|\EE_2 \left[\BB_{C,n,1}^M(\mathbf{u}) \BB_{C,n,2}^M(\mathbf{v}) \right] \right|
		& \le  \frac{1}{\sqrt{ \lfloor \lambda n \rfloor (n-\lfloor \lambda n \rfloor)} } \sum_{i=1}^{\lfloor \lambda n \rfloor} \sum_{j=\lfloor \lambda n \rfloor+1}^n | \Cov( \xi_{n,i},\xi_{n,j}) | \\
		& \le const \frac{1}{n}\sum_{i=1}^{c\cdot l(n)-1} i \left|v\left(\frac{i}{l(n)} \right)\right| =o(1).
\end{align*}
 Therefore, after some similar calculations as in the proof of Theorem \ref{theorem:asymptoticsS}, $\sigma_n^2=\EE_2[(\sum_{j=1}^n \zeta_{n,j})^2] \to \Var(c_1\BB_{C,1}(\mathbf{u}) + c_2 \BB_{C,2}(\mathbf{v}))$, $\PP_1$-almost surely. Since $| \zeta_{n,j} | =O(|\xi_{n,j}-1|\cdot n^{-1/2})$ uniformly in $\mathbf{u},\mathbf{v}$ and $\omega_1$, condition (2.1) in \cite{peligrad1996} is satisfied $\PP_1$-almost surely. Finally, we check the Lindeberg condition (2.2) in \cite{peligrad1996}, i.e., for all $\eps>0$,
\[
	\frac{1}{\sigma_n^2} \sum_{j=1}^n \EE_2[ \zeta_{n,j}^2 \mathbf{1}_{\{ |\zeta_{n,j}/\sigma_n| > \eps\} } ] \to 0 \quad  (\PP_1\text{-a.s.})
\]
for $n\to \infty$. The expression on the left-hand side of this display can be estimated uniformly in $\mathbf{u},\mathbf{v}$ and $\omega_1$ by
\[
	const\cdot \EE_2[ (\xi_{n,1}-1)^2 \mathbf{1}_{\{ |\xi_{n,1}-1| > const \sqrt{n} \}}] \le const \cdot \frac{1}{n} \EE_2[ (\xi_{n,1}-1)^4 ],
\]
which yields the assertion by the assumption on the finite centered moments of $\xi_{n,j}$. To conclude,
\[
	(\BB_{C,n,1}^M, \BB_{C,n,2}^M)  \weak (\BB_{C,1}, \BB_{C,2}) \quad  (\PP_1\text{-a.s.})
\]
and analogously to the proof of Theorem \ref{theorem:taperedmultiplier} this implies
\[
	(\widehat{\BB}_{C,n,1}^M, \widehat{\BB}_{C,n,2}^M)  \weak (\BB_{C,1}, \BB_{C,2}) \quad  (\PP_1\text{-a.s.}).
\]
Thus, by the continuous mapping Theorem
\begin{multline*}
	\sqrt{1-\lambda} \bigg\{ \widehat{\BB}_{C,n,1}^M - \sum_{i=1}^d D_iC(\mathbf{u}) \widehat{\BB}_{C,n,1}^M(\mathbf{u}^{(i)}) \bigg\}
		- \sqrt{\lambda} \bigg\{ \widehat{\BB}_{C,n,1}^M - \sum_{i=1}^d D_iC(\mathbf{u}) \widehat{\BB}_{C,n,1}^M(\mathbf{u}^{(i)}) \bigg\} \\
			\weak \sqrt{1-\lambda} \GG_{C,1}(\mathbf{u}) - \sqrt{\lambda} \GG_{C,2}(\mathbf{u}) \quad  (\PP_1\text{-a.s.}).
\end{multline*}
Again, the expression on the left-hand side is a ball-measurable element of $D([0,1]^d)$ which allows to translate the latter weak convergence into the Hoffmann-J\o rgensen conditional weak convergence. Finally, replace the true partial derivatives by their estimators and note that the difference is uniformly $o_P(1)$, which yields the assertion after an application of Lemma B.1 in \cite{buedetvol2011}.
\end{proof}

\begin{proof}[{\it Proof of Theorem \ref{theorem:kieferasymptotics}}]
Without loss of generality we may assume $\mathbf{X}_j=\mathbf{U}_j$. Under the conditions on the mixing rate, it follows from Theorem 2 in \citet{philpinzur80} that the process
\[
	\BB_{C,n}(\zeta,\mathbf{u}):=\frac{1}{\sqrt{n}} \sum_{j=1}^{\lfloor \zeta n \rfloor} \left\{ \mathbf{1}_{\{\mathbf{X}_j \leq \mathbf{u}\}} -C(\mathbf{u}) \right\}
\]
weakly converges in $\ell^\infty([0,1]^{d+1})$ to a $C$-Kiefer process $\BB_{C}(\zeta,\mathbf{u})$.
Hence, Condition~3.1 in~\cite{buecvolg2011} is satisfied and an application of Corollary 3.3 in that reference yields that
\[
	\tilde \Sb_{n}(\zeta,\mathbf{u}) = \frac{1}{\sqrt{n}} \sum_{j=1}^{\lfloor \zeta n \rfloor} \left\{ \mathbf{1}_{\{X_{j1} \le F_{n1}^{-1}(u_1), \dots, X_{jd} \le F_{nd}^{-1}(u_d) \}} - C_n(C_{n1}^{-1}(u_1), \dots, C_{nd}^{-1}(u_d)) \right\}
\]
weakly converges to $\Sb_C(\zeta,\mathbf{u})= \BB_C(\zeta,\mathbf{u})- \zeta \BB_C(1,\mathbf{u})$ in $\ell^\infty([0,1]^{d+1})$. Under the null hypothesis we can write
\[
	\Sb_n(\zeta,\mathbf{u}) = \frac{1}{\sqrt{n}} \sum_{j=1}^{\lfloor \zeta n \rfloor} \left\{ \mathbf{1}_{\{ \widehat{\mathbf{U}}_j \le \mathbf{u} \}} - \widehat{C}_n(\mathbf{u}) \right\}
\]
and the assertion follows from $\| \tilde \Sb_n - \Sb_n \|_\infty=O_P(n^{-1/2})$.
\end{proof}

\begin{proof}[{\it Proof of Proposition \ref{proposition:estimationSunspec}}]
Without loss of generality we may assume $\mathbf{X}_j=\mathbf{U}_j$. Suppose the null hypothesis holds and define
\[
	\BB_{C,n}^M(\zeta,\mathbf{u})
		=  \frac{1}{\sqrt{n}} \sum_{j=1}^{\lfloor \zeta n \rfloor} \left( \frac{\xi_{n,j}}{\bar \xi_{\lfloor\zeta n \rfloor}} -1 \right) \mathbf{1}_{\{\mathbf{X}_j \leq \mathbf{u}\}}
			=\frac{1}{\sqrt{n}} \sum_{j=1}^{\lfloor \zeta n \rfloor} \left( \frac{\xi_{n,j}}{\bar \xi_{\lfloor\zeta n \rfloor}} -1 \right) \left\{ \mathbf{1}_{\{\mathbf{X}_j \leq \mathbf{u}\}} - C(u) \right\}.
\]
We may replace $\bar \xi_{\lfloor\zeta n \rfloor}$ by $1$.
Similar as in the proofs of Proposition \ref{prop:estimationS} and Theorem \ref{theorem:taperedmultiplier} we begin by showing that
$\BB_{C,n}^M \weak \BB_C$, $\PP_1$- almost surely, where $\weak$ denotes weak convergence with respect to the multipliers. Convergence of the finite dimensional distributions follows along similar lines as in the previous proofs; the details are omitted for the sake of brevity. Tightness follows by tightness of $\BB_{C,n}^M(1,\cdot)$ and by analogous arguments as in the proof of Theorem 2.12.1 in \cite{vaartwellner96}.

Thus, as in the proof of Theorem \ref{theorem:taperedmultiplier}, $\widehat{\BB}_{C,n}^M \weakcondp{\xi} \BB_C$; and hence  $\widehat{\Sb}_{n}^M \weakcondp{\xi} \Sb$ by the Lipschitz-continuous mapping Theorem for the bootstrap, see Proposition 10.7 in \cite{kosorok08}.
Under the alternative, again replacing $\bar \xi_{\lfloor\zeta n \rfloor}$ by $1$, we have
\[
	\frac{1}{n}\EE_2 \left( \sum_{j=1}^{\lfloor \zeta n \rfloor} \left( {\xi_{n,j}} -1 \right) \mathbf{1}_{\{\mathbf{X}_j \leq \mathbf{u}\}} \right)^2 \le \frac{l_n}{n} \sum_{k=1}^{l_n} (1-k/l_n) v(k/l_n) =o(1)
\]
by dominated convergence, which yields $\|\widehat{\BB}_{C,n}^M(\zeta,\mathbf{u})\|_\infty = o_{\PP_2}(1)$ uniformly in $\zeta, \mathbf{u}$ and $\omega_1$.

\end{proof}



\bibliographystyle{apalike}
\bibliography{bibliography_jmva}

\end{document}